\pdfoutput=1
\documentclass[11pt]{amsart}
\usepackage[margin=1.1in]{geometry}
\usepackage{amsmath,amssymb,amsthm,graphicx,booktabs}

\usepackage[dvipsnames]{xcolor}
\usepackage[section]{placeins}
\usepackage{hyperref}
\hypersetup{colorlinks=true,linkcolor=blue!50!black,citecolor=blue!50!black,urlcolor=blue!50!black}

\newtheorem{theorem}{Theorem}
\newtheorem{proposition}[theorem]{Proposition}
\newtheorem{lemma}[theorem]{Lemma}
\newtheorem{corollary}[theorem]{Corollary}
\newtheorem{conjecture}[theorem]{Conjecture}
\newtheorem{mainthm}{Theorem}

\theoremstyle{remark}
\newtheorem{remark}[theorem]{Remark}

\renewcommand{\Re}{\operatorname{Re}}
\renewcommand{\Im}{\operatorname{Im}}

\title[Gram's law for moving targets]{Gram's law for moving targets of the Riemann zeta-function}
\author{Ryan Goulden}
\address{Independent researcher, Crown Point, Indiana, USA}
\email{rgoulden000@gmail.com}
\date{}

\begin{document}

\begin{abstract}
We study the solutions of $\zeta(s)=f(s)$ for targets $f$ that grow with $|s|$, and Gram's law for them: the statement that each shifted Gram interval contains exactly one solution. The deviation from Gram's law is an integer phase-slip counter, equal up to rounding to $-\tfrac12S(t;x_c)$, where $S(t;x)=\pi^{-1}\arg\zeta(x+it)$ and $x_c\approx\tfrac12+\log|f(it)|/\log t$ is a balance abscissa determined by the target.

If $|f|\asymp t^{\kappa}$ with $\kappa>\kappa_1=0.5339\ldots$, where $\sum_p\arcsin(p^{-1/2-\kappa_1})=\pi$, and $\arg f$ converges, then every sufficiently high shifted Gram interval contains exactly one solution, which is simple; this localises the fixed points of $\zeta$. The threshold is sharp: for $f=c(-s)^{\kappa}$ with $\tfrac12<\kappa<\kappa_1$, Gram's law fails in $\asymp T$ intervals up to height $T$, and for $0<\kappa\le\tfrac12$ it fails infinitely often but with density zero. In both ranges the counter has an explicit limiting distribution.

For targets growing more slowly than any power, such as $\log(-s)$, the solutions lie on average $(1+o(1))\log|f(iT)|/\log T$ to the left of the critical line. If $\log|f(iT)|/\sqrt{\log\log T}\to\infty$, the proportion of solutions to the right of the line tends to $0$, whereas for fixed targets Selberg conjectured a quarter, and the counting error at the shifted Gram points is asymptotically Gaussian with variance $(8\pi^2)^{-1}\log(\log T/\log|f(iT)|)+O(1)$. For $\log(-s)$ the solutions follow neither $\theta(T)/\pi$ nor the zeros of $\zeta$. All results are unconditional; numerical computations for $\log(-s)$ illustrate how slowly they take effect.
\end{abstract}

\maketitle

\section{Introduction}

For a fixed $a\in\mathbb C$, the solutions of $\zeta(s)=a$ (the \emph{$a$-points}) have been studied since Landau. They are as numerous as the zeros, and almost all of them cluster around the critical line \cite{Levinson}. For fixed $a\neq0$, Selberg conjectured, on heuristic grounds, that about three quarters of the nontrivial $a$-points lie to the left of $\sigma=\tfrac12$ and about one quarter to the right. On RH he showed, in unpublished work, that about half of them lie to the left at distance of order $\sqrt{\log\log T}/\log T$ (see \cite{Tsang,Lester}).

Equations $\zeta(s)=f(s)$ with a non-constant target $f$ (\emph{moving targets}, as in the Nevanlinna theory of small functions \cite{Yamanoi}) have been treated with Nevanlinna methods. Li and Steuding \cite{LiSteuding} showed that the Nevanlinna counting function of the solutions of $\zeta^{(m)}(s)=P(s)$, for a polynomial $P$, satisfies $N(r)=\frac1\pi r\log r+O(r)$; see also \cite{LiSteudingSuzuki}. Cardwell and Ye \cite{CardwellYe} used Littlewood's lemma to bound sums of $\beta-b$ over the solutions of $L(s)=R(s)$, for rational targets $R$ and a class of $L$-functions containing the Selberg class, with error $O(T)$. In the other direction, $a$-points of the factor $\Delta=\chi$ in the functional equation, and their relation to $a$-points of $\zeta$, have been studied by Kalpokas--Steuding \cite{KalpokasSteuding}, Steuding--Suriajaya \cite{SteudingSuriajaya} and Sourmelidis--Steuding--Suriajaya \cite{SSS}.

Moving targets also make the failure of Gram's law adjustable. For the zeros of $\zeta$, Gram's law fails infinitely often because $S(t)$ is unbounded, and in fact for a positive proportion of Gram intervals \cite{Trudgian}. For a growing target the same mechanism operates, with $S(t)$ replaced by the argument of $\zeta$ at an abscissa $x_c>\tfrac12$ that the target selects (Section~\ref{sec:balance}): the faster $f$ grows, the further to the right $x_c$ lies, and the better behaved $\arg\zeta$ is there. Varying the target therefore moves the problem from the half-plane of absolute convergence to the critical line, and it exposes two transitions: at the explicit exponent $\kappa_1=0.5339\ldots$, and where $\log|f|$ becomes comparable with $\sqrt{\log\log t}$.

\subsection{An observation}
This paper began with a numerical observation about $\zeta(s)=\log(-s)$ (principal branch). With $-1.5\le\sigma\le3$ and $1\le t\le400$ there are exactly $202$ solutions, the same as the number of zeros of $\zeta$, and all but $8$ of them lie to the left of the critical line. The explanation is general:
\begin{enumerate}
\item a target with $|f|\gg1$ forces solutions to the left of $\sigma=\tfrac12$, where the factor $|\chi(s)|=(t/2\pi)^{1/2-\sigma}$ supplies the size;
\item the mirror point $1-\rho$ then lies to the right, where $\arg\zeta$ is small; and
\item the functional equation turns this into phase locking with shifted Gram points.
\end{enumerate}
Taking logarithms in $f(\rho)=\chi(\rho)\zeta(1-\rho)$ makes this precise (Theorem~\ref{thm:phase}). A solution $\rho=\beta+i\gamma$ satisfies
\begin{align*}
(\tfrac12-\beta)\log\frac{\gamma}{2\pi}&=\log|f(\rho)|-\log|\zeta(1-\rho)|+O(\gamma^{-1}),\\
2\theta(\gamma)&\equiv\arg\zeta(1-\rho)-\arg f(\rho)+O(\gamma^{-1})\pmod{2\pi}.
\end{align*}
The first relation places the solutions near the line $\sigma=\tfrac12-\log|f|/\log(\gamma/2\pi)$. The second ties their ordinates to the Gram phase $2\theta(\gamma)$, up to the phase of $\zeta$ at the mirror point.

\subsection{The balance-point principle}\label{sec:balance}
\begin{figure}[t]
\centering
\includegraphics[width=0.9\linewidth]{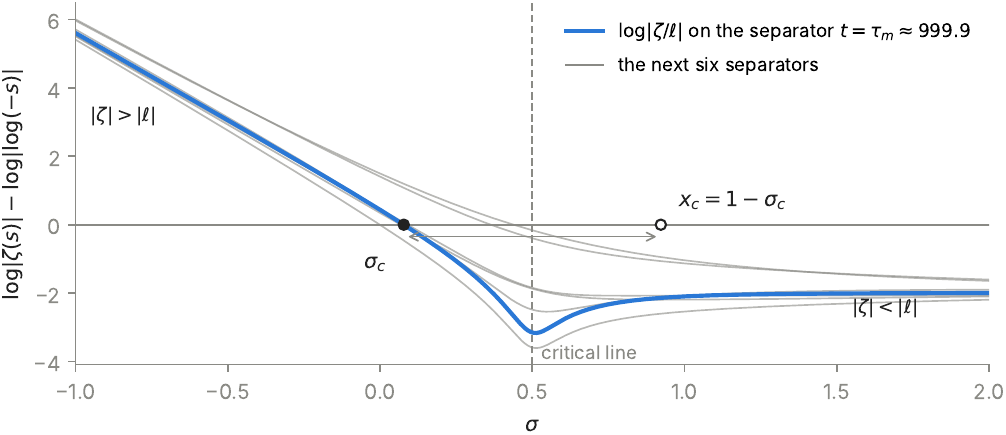}
\caption{The balance point for $f(s)=\log(-s)$. Along the separator $t=\tau_m\approx999.9$, $\log|\zeta/f|$ decreases through $0$ exactly once, at $\sigma_c\approx0.08$; the balance abscissa is the mirror point $x_c=1-\sigma_c\approx0.92$. The thin curves show the next six separators, each with a single crossing.}\label{fig:balance}
\end{figure}

Fix a large height $\tau$. On the line $t=\tau$ the ratio $R(\sigma)=|\zeta(\sigma+i\tau)/f(\sigma+i\tau)|$ is large on the left, where $|\chi|$ dominates, and small on the right. For the targets studied here it crosses $1$ exactly once (at almost all heights, when the crossing is close to the critical line), at an abscissa $\sigma_c$. We call its mirror image
\[
x_c=1-\sigma_c=\tfrac12+\frac{\log|f(i\tau)|-\log|\zeta(x_c+i\tau)|}{\log(\tau/2\pi)}+O\Bigl(\frac1\tau\Bigr)
\]
the \emph{balance abscissa}. Let $\varphi_0=\lim\arg f$, define the shifted Gram points $\tau_m$ by $\theta(\tau_m)=(m+\tfrac12)\pi-\tfrac12\varphi_0$, and let $n_m$ be the number of solutions with $\tau_{m-1}<t\le\tau_m$. Gram's law is the statement $n_m=1$. We say that the solutions \emph{lock} onto the shifted Gram points if $n_m=1$ for all large $m$. Our structural result (Propositions~\ref{prop:counter} and~\ref{prop:ident}) is that
\[
N_f(\tau_m)=m+C_f-k_m,\qquad k_m=\text{the integer nearest to }-\tfrac12S(\tau_m;x_c)+o(1),
\]
where $N_f$ counts the solutions and $S(t;x)=\pi^{-1}\arg\zeta(x+it)$. Thus $n_m=1+k_{m-1}-k_m$: Gram's law holds exactly as long as the integer $k_m$ is constant, and it fails when $S(\cdot\,;x_c)$ crosses an odd integer, which we call a \emph{phase slip}. Figure~\ref{fig:balance} shows the crossing along a separator for $f(s)=\log(-s)$.

For $x>1$ the argument of $\zeta(x+it)$ is bounded by $A(x)=\sum_p\arcsin(p^{-x})$, and $A(x)<\pi$ exactly when $x>\sigma_1=1.0339\ldots$. Put $\kappa_1=\sigma_1-\tfrac12$. As $x_c$ moves to the left, the phase at $x_c$ needs successively deeper tools: the bound $A(x_c)<\pi$, then Kronecker's theorem for $1<x_c<\sigma_1$, the Bohr--Jessen limit theorem for fixed $x_c\in(\tfrac12,1]$, and the Selberg--Tsang central limit theorem when $x_c\to\tfrac12$. This gives the following picture.
\begin{center}
\small
\begin{tabular}{@{}llll@{}}
\toprule
target & balance abscissa & counter $k_m$ & Gram's law\\
\midrule
$|f|\asymp t^{\kappa}$, $\kappa>\kappa_1$ & $x_c>\sigma_1$ & eventually constant & all large $m$\\
$c(-s)^{\kappa}$, $\tfrac12<\kappa<\kappa_1$ & $1<x_c<\sigma_1$ & limit law, bounded support & fails $\asymp T$ times\\
$c(-s)^{\kappa}$, $0<\kappa\le\tfrac12$ & $\tfrac12<x_c\le1$ & limit law on $\mathbb Z$ & fails i.o., density $0$\\
slowly growing${}^{*}$ & $x_c\to\tfrac12$ & Gaussian, variance $\to\infty$ & fails i.o.\\
$a$ fixed & none & --- & no locking\\
\bottomrule
\multicolumn{4}{@{}l}{\rule{0pt}{2.6ex}\footnotesize${}^{*}$\,$\log|f|=o(\log t)$ and $\log|f|/\sqrt{\log\log t}\to\infty$, for instance $f=\log(-s)$.}
\end{tabular}
\end{center}
The first relation above also governs the horizontal distribution. Power targets push the solutions to $\sigma\approx\tfrac12-\kappa$. Slowly growing targets such as $\log(-s)$ push them to the left by about $\log|f|/\log t$. Once $\log|f|$ is large compared with $\sqrt{\log\log t}$, the typical size of $\log|\zeta|$ near the critical line, this leaves almost no solutions to the right of the line (Theorem~\ref{thm:right}), so there is no analogue of Selberg's conjectured three-quarters/one-quarter split.

\subsection{Main results}\label{sec:main}
We write $Y_x=-\sum_p\arg(1-p^{-x}e^{i\vartheta_p})$ for $x>\tfrac12$, where the $\vartheta_p$ are independent and uniformly distributed on $[0,2\pi)$; the series converges almost surely. Its law is the limiting distribution of $\arg\zeta(x+it)$.

\begin{mainthm}[Locking; Theorem~\ref{thm:gram}, Corollaries~\ref{cor:count} and~\ref{cor:poly}]
Let $\kappa>\kappa_1$, and let $f$ satisfy the hypotheses \textup{(H1)}, \textup{(H2)} of Section~\ref{sec:big}, with $\arg f\to\varphi_0$. Then for all large $m$ the strip $\tau_{m-1}<t\le\tau_m$ contains exactly one solution of $\zeta(s)=f(s)$, and it is simple. Hence $N_f(T)=(\theta(T)+\varphi_0/2)/\pi+O(1)=N(T)-S(T)+O(1)$. For a polynomial $P$ of degree $d\ge1$, the solutions of $\zeta(s)=P(s)$ satisfy $\beta=\tfrac12-d+O(1/\log\gamma)$. In particular the non-real fixed points of $\zeta$ have real parts tending to $-\tfrac12$, and there are $\frac{T}{2\pi}\log\frac{T}{2\pi e}+O(1)$ of them with $0<t\le T$.
\end{mainthm}

\begin{mainthm}[Sharpness; Theorems~\ref{thm:sharp} and~\ref{thm:order}]
Let $\tfrac12<\kappa<\kappa_1$ and $x_0=\tfrac12+\kappa$. For $f(s)=c(-s)^{\kappa}$ with $c\ne0$,
\[
\sum_{\tau_m\le T}|n_m-1|\asymp T,\qquad \frac1M\#\{m\le M:\ k_m=j\}\to\mathbb P\bigl((2j-1)\pi<-Y_{x_0}<(2j+1)\pi\bigr)\quad(j\in\mathbb Z).
\]
The number of failing intervals up to height $T$ is $\gg T$ and $o(T\log T)$ for every target satisfying \textup{(H1$'$)}, \textup{(H2)}.
\end{mainthm}

\begin{mainthm}[Small exponents; Theorem~\ref{thm:smallkappa}]
Let $f(s)=c(-s)^{\kappa}$ with $c\ne0$ and $0<\kappa\le\tfrac12$. Then $k_m$ has the limiting distribution of Theorem~B, now with $x_0=\tfrac12+\kappa\in(\tfrac12,1]$, and every integer has positive limiting probability. Consequently Gram's law fails for infinitely many $m$, although the number of failures with $\tau_m\le T$ is $o(T\log T)$.
\end{mainthm}

\begin{mainthm}[The logarithmic target]
Let $N_H(T)$ count the solutions of $\zeta(s)=\log(-s)$ with $-1<\sigma<2$ and $0<t\le T$, and put $S_H(T)=N_H(T)-\theta(T)/\pi$.
\begin{enumerate}
\item[(i)] $N_H(T)=\theta(T)/\pi+O(\log T)$, and the sum of $\beta-\tfrac12$ over the solutions with $0<\gamma\le T$ is $-\frac{T}{2\pi}\log\log T+O(T/\log T)$ \textup{(}Theorems~\ref{thm:countH} and~\ref{thm:drift}\textup{)}.
\item[(ii)] The proportion of solutions with $T<\gamma\le2T$ and $\beta>\tfrac12$ is $\ll1/\log\log T$ \textup{(}Theorem~\ref{thm:right}\textup{)}.
\item[(iii)] $S_H(\tau_m)/\sigma_T$ is asymptotically standard Gaussian over the shifted Gram points, where $\sigma_T^2\sim\log\log T/8\pi^2$. In particular Gram's law fails for infinitely many $m$ \textup{(}Theorem~\ref{thm:gramfail}\textup{)}.
\item[(iv)] $|N_H(\tau_m)-N(\tau_m)|\ge c\sqrt{\log\log T}$ for a positive proportion of $m$. With (iii), the solutions follow neither $\theta(T)/\pi$ nor the zeros of $\zeta$ \textup{(}Corollary~\ref{cor:neither}\textup{)}.
\end{enumerate}
\end{mainthm}

\begin{mainthm}[Slowly growing targets; Theorems~\ref{thm:countH}--\ref{thm:right} and~\ref{thm:gramfail}]
Let $f$ satisfy the hypotheses of Section~\ref{sec:slow}, and put $L_f(T)=\log|f(iT)|$. Then $N_f(T)=\theta(T)/\pi+O(\log T)$, the mean of $\beta-\tfrac12$ over the solutions with $T<\gamma\le2T$ is $-(1+o(1))L_f(T)/\log T$, and
\[
\#\{T<\gamma\le2T:\ \beta>\tfrac12\}\ll\frac{T\log T\,\log\log T}{L_f(T)^2}+\log^2T .
\]
If moreover $\arg f(\sigma+it)\to\varphi_0$ uniformly for $-1\le\sigma\le2$, and $L_f(T)/\sqrt{\log\log T}\to\infty$, then the counting error $S_f(\tau_m)=N_f(\tau_m)-(\theta(\tau_m)+\tfrac12\varphi_0)/\pi$ at the shifted Gram points is asymptotically Gaussian, with variance $(8\pi^2)^{-1}\log(\log T/L_f(T))+O(1)$. In particular Gram's law fails for infinitely many $m$.
\end{mainthm}

The variance in Theorem~E is $\sim\log\log T/8\pi^2$ whenever $\log L_f(T)=o(\log\log T)$, as for $c\,\ell(s)^A$, and it is $\sim(1-\nu)\log\log T/8\pi^2$ when $L_f(T)\asymp(\log T)^\nu$ with $0<\nu<1$, as for $\exp(\ell(s)^\nu)$. The hypothesis $L_f(T)/\sqrt{\log\log T}\to\infty$ compares the push to the left with the typical size of $\log|\zeta|$ near the critical line. Below it we expect a positive proportion of the solutions to lie to the right of the line (Remark~\ref{rem:split}).

\subsection{Relation to earlier work}
The phase identity \eqref{eq:phase} below is the logarithmic form of the relation $\chi(\rho)=f(\rho)/\zeta(1-\rho)$. For constant $f=a$ this relation underlies \cite{SSS}. The analogous congruence for the $a$-points of $\Delta=\chi$ is \cite[eq.~(5)]{SteudingSuriajaya}. We claim no novelty for it. Theorems~\ref{thm:sharp}--\ref{thm:order} are moving-target analogues of the failure of Gram's law for the zeros \cite{Trudgian}, with $S(t)$ replaced by an integer phase-slip counter. As far as we can determine, the following are new: the per-interval locking of Theorem~A and the sharpness of its threshold; the balance-point description of the counter and the limiting distributions in Theorems~B and~C; and, for slowly growing targets, the drift law, the result on solutions to the right of the critical line and the central limit theorem for the counting error, including the analysis of Gram's law for $\log(-s)$ (Theorems~D and~E). The Riemann--von Mangoldt formula for slowly growing targets follows by the standard argument. For polynomial targets, the real parts of the solutions follow at once from \eqref{eq:modulus}; what is new is their localisation to individual shifted Gram intervals and the bounded counting error. The main tool for the results on slowly growing targets is a bound for the maximum of $\log|\zeta|$ over the right half of a short horizontal strip (Lemma~\ref{lem:max}), which we deduce from a Poisson majorant for $\log|\zeta|$, the Hardy--Littlewood maximal inequality and the moment bounds of Heap, Radziwi\l\l\ and Soundararajan. The results of \cite{LiSteudingSuzuki} concern the Nevanlinna characteristic of $\zeta$, with error $O(\log r/r)$, and the number of fixed points in discs; they do not localise individual solutions. Regularity results of a different kind, with simple zeros of $\zeta^{(k)}$ confined to explicit regions, are known for the derivatives of $\zeta$ \cite{BinderPauliSaidak}, and uniqueness questions for $L$-functions with moving targets are studied in \cite{WangWangYao}.

\subsection{Organisation}
Section~\ref{sec:prelim} collects preliminaries, and Section~\ref{sec:counter} sets up the phase-slip counter for a general target. Sections~\ref{sec:big} and~\ref{sec:smallkappa} treat power-like targets with $\kappa>\tfrac12$ and $0<\kappa\le\tfrac12$, Sections~\ref{sec:slow} and~\ref{sec:gramlog} slowly growing targets, Section~\ref{sec:numerics} the computations, and Section~\ref{sec:open} open problems.

\subsection*{Notation}
$s=\sigma+it$, with $t>0$ large. $\chi(s)=2^s\pi^{s-1}\sin(\pi s/2)\Gamma(1-s)$, so that $\zeta(s)=\chi(s)\zeta(1-s)$. $\theta$ is the Riemann--Siegel theta function, with $\chi(\tfrac12+it)=e^{-2i\theta(t)}$. $N(T)=\theta(T)/\pi+1+S(T)$. $\log$ is the principal branch. The letter $p$ always denotes a prime, and sums and products over $p$ run over all primes. For $\Re w>1$, $\arg\zeta(w)=\Im\log\zeta(w)=\Im\sum_{n}\Lambda(n)n^{-w}/\log n$. For $x\ge\tfrac12$ we put $S(t;x)=\pi^{-1}\arg\zeta(x+it)$, with $\arg$ defined by continuous variation from $2+it$ (when no zero lies on the segment), so that $S(t;\tfrac12)=S(t)$. Write $\ell(s)=\log(-s)$. For a finite set $A$, $\#A$ denotes the number of its elements, and $\operatorname{meas}$ denotes Lebesgue measure.

\section{Preliminaries}\label{sec:prelim}

\subsection{The factor \texorpdfstring{$\chi$}{chi} and the phase identity}

\begin{lemma}\label{lem:chi}
Uniformly for $|\sigma|\le A$ and $t\ge2$,
\[
\log\chi(\sigma+it)=(\tfrac12-\sigma)\log\frac{t}{2\pi}-2i\theta(t)+O_A(t^{-1}).
\]
\end{lemma}
\begin{proof}
By Stirling's formula for $\Gamma'/\Gamma$, $\chi'/\chi(s)=-\log(t/2\pi)+O(1/t)$ uniformly in the strip. Integrate in $\sigma$ from $\tfrac12$ and use $\chi(\tfrac12+it)=e^{-2i\theta(t)}$.
\end{proof}

\begin{lemma}\label{lem:left}
For $\sigma\le0$ and $t\ge1$,
\[
\log|\chi(\sigma+it)|\ge(\tfrac12-\sigma)\log\frac{|1-s|}{2\pi e}-C
\]
for an absolute constant $C$. Consequently, suppose $|f(s)|\le C_1|s|^{K}$ for $t\ge t_0$, and $\sigma_L<\min(0,\tfrac12-K)$. Then $|\zeta(s)|>|f(s)|$ for all $\sigma\le\sigma_L$ and $t\ge t_1(f,\sigma_L)$.
\end{lemma}
\begin{proof}
We have $|\sin(\pi s/2)|\ge\tfrac14e^{\pi t/2}$ for $t\ge1$. Put $w=1-s$, so $\Re w\ge1$ and $\arg w\in(-\pi/2,0)$. Stirling's formula gives
\[
\log|\Gamma(w)|=(\Re w-\tfrac12)\log|w|-\Im w\cdot\arg w-\Re w+O(1),
\]
and here $-\Im w\cdot\arg w=t\arg w\ge-\pi t/2$. Adding $\sigma\log2\pi-\log\pi$ gives the first claim.

For the second, note that $|\zeta(1-s)|\ge\zeta(2-2\sigma)/\zeta(1-\sigma)\ge c(\sigma_L)>0$ for $\sigma\le\sigma_L$. Put $B=\tfrac12-\sigma_L>K$. Since $\tfrac12-\sigma\ge B$ and $\log(|1-s|/2\pi e)>0$, we get
\[
\log|\zeta(s)|-\log|f(s)|\ge(B-K)\log|s|-O_{\sigma_L}(1),
\]
which is positive for $t$ large.
\end{proof}

\begin{theorem}[Phase identity]\label{thm:phase}
Let $f$ be analytic near $\rho=\beta+i\gamma$, with $|\beta|\le A$, $\gamma\ge2$, and $\zeta(\rho)=f(\rho)\ne0$. Then
\begin{align}
2\theta(\gamma)&\equiv \arg\zeta(1-\rho)-\arg f(\rho)+O_A(\gamma^{-1})\pmod{2\pi},\label{eq:phase}\\
(\tfrac12-\beta)\log\frac{\gamma}{2\pi}&=\log|f(\rho)|-\log|\zeta(1-\rho)|+O_A(\gamma^{-1}).\label{eq:modulus}
\end{align}
\end{theorem}
\begin{proof}
Take logarithms in $f(\rho)=\chi(\rho)\zeta(1-\rho)$ and apply Lemma~\ref{lem:chi}.
\end{proof}

\subsection{The function \texorpdfstring{$A(\sigma)$}{A(sigma)}}

For $\sigma>1$ put
\[
A(\sigma)=\sum_p\arcsin\bigl(p^{-\sigma}\bigr).
\]
This function is continuous and strictly decreasing, with $A(\sigma)\to\infty$ as $\sigma\to1^+$. Let $\sigma_1=1.0339080723\ldots$ be the solution of $A(\sigma_1)=\pi$, and put $\kappa_1=\sigma_1-\tfrac12=0.5339080723\ldots$. (Numerically, $A(\sigma)=\sum_{j\ge0}\binom{2j}{j}4^{-j}(2j+1)^{-1}\mathcal P((2j+1)\sigma)$, where $\mathcal P(s)=\sum_pp^{-s}$ is the prime zeta function.)

\begin{lemma}\label{lem:arcsin}
Let $\sigma>1$.
\begin{enumerate}
\item[(a)] $|\arg\zeta(\sigma+it)|\le A(\sigma)$ for all $t$.
\item[(b)] Let $\sigma_0>1$, $|v|<A(\sigma_0)$ and $\delta>0$. Then there are $h,d,L>0$ such that the set
\[
\bigl\{t:\ |\arg\zeta(\sigma+it)-v|<\delta\ \text{for all }\sigma\in[\sigma_0-h,\sigma_0+h]\bigr\}
\]
contains an interval of length $d$ inside every interval of length $L$.
\end{enumerate}
\end{lemma}
\begin{proof}
(a) For $|z|=r<1$ we have $|\arg(1-z)|\le\arcsin r$. Since $\arg\zeta(w)=-\sum_p\arg(1-p^{-w})$ for $\Re w>1$, the bound follows by summing.

(b) Choose $y$ and $h$ with $\sum_{p>y}\arcsin(p^{-\sigma_0+h})<\delta/4$. Then $v$ lies in the interior of the range of the continuous map
\[
(\vartheta_p)_{p\le y}\mapsto-\sum_{p\le y}\arg\bigl(1-p^{-\sigma_0}e^{i\vartheta_p}\bigr)
\]
from the torus, whose range is $[-A_y,A_y]$, where $A_y=\sum_{p\le y}\arcsin(p^{-\sigma_0})\to A(\sigma_0)$. So we may fix phases $\vartheta_p$ with value $v$. By uniform continuity, after shrinking $h$ there is $\nu>0$ such that $|t\log p+\vartheta_p|<\nu\pmod{2\pi}$ for all $p\le y$ implies the conclusion with $\delta$ replaced by $3\delta/4$. The $\log p$ are linearly independent over $\mathbb Q$, so by Kronecker's theorem the set of $t$ satisfying these simultaneous conditions with $\nu/2$ is relatively dense. Each such $t$ lies in an interval of length $d=\nu/(2\log y)$ on which the conditions hold with $\nu$.
\end{proof}

\subsection{A winding lemma}

When $\arg\zeta$ is not bounded by $\pi$, the winding along a separator is no longer determined by the endpoints. The following elementary lemma computes it exactly whenever the modulus ratio crosses $1$ only once. Throughout, $\operatorname{Arg}$ denotes the principal argument; recall that $|\operatorname{Arg}(1+z)|\le\arcsin|z|$ for $|z|<1$.

\begin{lemma}[Single crossing]\label{lem:cross}
Let $a<b$, and let $u\colon[a,b]\to\mathbb C$ be continuous with $1+u\ne0$ on $[a,b]$. Suppose there is $\sigma^*\in(a,b)$ with $|u|\ge1$ on $[a,\sigma^*]$ and $|u|\le1$ on $[\sigma^*,b]$, and let $\psi$ be a continuous argument of $u$ on $[a,\sigma^*]$. Then $\psi(\sigma^*)\notin\pi+2\pi\mathbb Z$, and the change of $\arg(1+u)$ along $[a,b]$ is
\[
-\psi(a)+2\pi j+\operatorname{Arg}\bigl(1+u(b)\bigr)-\operatorname{Arg}\bigl(1+u(a)^{-1}\bigr),
\]
where $j$ is the integer nearest to $\psi(\sigma^*)/2\pi$.
\end{lemma}
\begin{proof}
If $|z|\le1$ and $z\ne-1$, then $\Re(1+z)\ge1-|z|\ge0$, with equality only for $z=-1$. So $\operatorname{Arg}(1+z)\in(-\tfrac\pi2,\tfrac\pi2)$ and it depends continuously on $z$. Let $\omega$ be a continuous argument of $1+u$ on $[a,b]$. On $[a,\sigma^*]$ we have $1+u=u(1+u^{-1})$ with $|u^{-1}|\le1$, so $\omega=\psi+\operatorname{Arg}(1+u^{-1})+2\pi n_1$ for a constant integer $n_1$. On $[\sigma^*,b]$ we have $\omega=\operatorname{Arg}(1+u)+2\pi n_2$ for a constant integer $n_2$.

At $\sigma^*$ both inequalities hold, so $|u(\sigma^*)|=1$ and $u(\sigma^*)=e^{i\psi(\sigma^*)}$. Since $u(\sigma^*)\ne-1$, we can write $\psi(\sigma^*)=2\pi j+\psi'$ with $|\psi'|<\pi$. Then $1+e^{\pm i\psi'}=2\cos(\psi'/2)\,e^{\pm i\psi'/2}$ with $\cos(\psi'/2)>0$. Hence $\operatorname{Arg}(1+u(\sigma^*)^{-1})=-\psi'/2$ and $\operatorname{Arg}(1+u(\sigma^*))=\psi'/2$. Comparing the two expressions for $\omega(\sigma^*)$ gives $n_2-n_1=j$, and then $\omega(b)-\omega(a)$ is as stated.
\end{proof}

\subsection{The size of \texorpdfstring{$\zeta$}{zeta} to the right of the critical line}\label{sec:rightsize}

For $y>0$ let $P_y(u)=y/(\pi(y^2+u^2))$ be the Poisson kernel of the half-plane, and write $(P_y*g)(t)=\int_{\mathbb R}P_y(t-v)g(v)\,dv$.

\begin{lemma}[Poisson majorant]\label{lem:poisson}
Let $0<a<1$ and $F(s)=(s-1)\zeta(s)/s$. For $\sigma>a$ and $t\in\mathbb R$,
\begin{equation}\label{eq:poisson}
\log|F(\sigma+it)|\le(P_{\sigma-a}*g)(t),\qquad g(v)=\log|F(a+iv)|.
\end{equation}
Consequently, uniformly for $\tfrac14\le a\le\tfrac34$, $a<\sigma\le3$ and $t\ge2$,
\begin{equation}\label{eq:poisson2}
\log|\zeta(\sigma+it)|\le(P_{\sigma-a}*G_a)(t)+O(1/t),\qquad G_a(v)=\log|\zeta(a+iv)|.
\end{equation}
\end{lemma}
\begin{proof}
The function $F$ is analytic in $\Re s>0$ and tends to $1$ as $\sigma\to\infty$, and the convexity bound for $\zeta$ gives $\log|F(s)|\le A\log(2+|s|)$ for $\Re s\ge a$, with $A$ depending on $a$. So $u=\log|F|$ is subharmonic in $\Re s>0$, and $g$ is continuous as a map into $[-\infty,\infty)$, with $g(v)\le A\log(2+|v|)$.

Fix $M>0$ and put $g_M=\max(g,-M)$. This is continuous, and $-M\le g_M(v)\le A\log(2+|v|)$. Hence $U_M(\sigma+it)=(P_{\sigma-a}*g_M)(t)$ is harmonic in $\sigma>a$, satisfies $-M\le U_M(s)\ll\log(2+|s|)$, and extends continuously to the line $\sigma=a$ with boundary values $g_M$. Put $w=u-U_M$. Then $w$ is subharmonic in $\sigma>a$, $w(s)\le A\log(2+|s|)+M$, and for every $v_0$
\[
\limsup_{s\to a+iv_0}w(s)\le u(a+iv_0)-g_M(v_0)\le0,
\]
because $u$ is upper semicontinuous and $u(a+iv_0)=g(v_0)\le g_M(v_0)$. For $\varepsilon>0$ let $w_\varepsilon(s)=w(s)-\varepsilon\Re(s-a)^{1/2}$, with the principal branch. Since $\Re(s-a)^{1/2}\ge|s-a|^{1/2}/\sqrt2\ge0$ in $\sigma\ge a$, the function $w_\varepsilon$ is subharmonic in $\sigma>a$, has boundary $\limsup\le0$ on the line $\sigma=a$, and tends to $-\infty$ as $|s|\to\infty$ in $\sigma\ge a$. By the maximum principle on the half-discs $\{\sigma>a,\ |s-a|<R\}$ with $R$ large, $w_\varepsilon\le0$. Letting $\varepsilon\to0$ gives $u\le U_M$. As $M\to\infty$, $g_M$ decreases to $g$ and $g_M\le g_1$ for $M\ge1$, with $P_{\sigma-a}*g_1$ finite. So $U_M$ decreases to $P_{\sigma-a}*g$ by monotone convergence, which proves \eqref{eq:poisson}.

For \eqref{eq:poisson2}, note that $\log|(s-1)/s|=O(1/t)$ for $t\ge2$. On the line $\sigma=a$, $\log|(a-1+iv)/(a+iv)|\le C\min(1,1/|v|)$ for $\tfrac14\le a\le\tfrac34$. Write $y=\sigma-a\le3$. The part of $(P_y*\min(1,1/|\cdot|))(t)$ with $|v|\ge t/2$ is at most $2/t$, and the part with $|v|<t/2$ is at most $t\cdot4y/(\pi t^2)$. So this convolution is $O(1/t)$, and \eqref{eq:poisson2} follows from \eqref{eq:poisson}.
\end{proof}

\begin{lemma}[Second moment]\label{lem:moment}
Uniformly for $\tfrac12\le\sigma\le3$,
\[
\int_T^{2T}\bigl(\log^+|\zeta(\sigma+it)|\bigr)^2\,dt\ll T\log\log T .
\]
\end{lemma}
\begin{proof}
First let $\sigma=\tfrac12$. By \cite[Thm.~1]{HRS}, $\int_T^{2T}|\zeta(\tfrac12+it)|^{2k}\,dt\ll T(\log T)^{k^2}$ uniformly for $(\log T)^{-1/2}\le k\le2$; this uniformity is recorded in \cite[Cor.~1]{HS}. Take $k=(\log\log T)^{-1/2}$, so that $(\log T)^{k^2}=e$. Since $\log y\le y^k/(ek)$ for $y\ge1$, we have $(\log^+y)^2\le(ek)^{-2}y^{2k}$ for $y>0$, and the claim follows.

Now let $\tfrac12<\sigma\le3$ and $y=\sigma-\tfrac12$. Write $G=G_{1/2}$. By \eqref{eq:poisson2} with $a=\tfrac12$, and since $P_y*G\le P_y*G^+$ and $x\mapsto x^+$ is increasing, $\log^+|\zeta(\sigma+it)|\le(P_y*G^+)(t)+C/t$. As $P_y$ is a probability density, the Cauchy--Schwarz inequality gives $(P_y*G^+)^2\le P_y*(G^+)^2$, so
\[
\int_T^{2T}\bigl(\log^+|\zeta(\sigma+it)|\bigr)^2\,dt\le2\int_{\mathbb R}G^+(v)^2\Bigl(\int_T^{2T}P_y(t-v)\,dt\Bigr)dv+O(1/T).
\]
The inner integral is at most $1$. For $v\notin[T/2,3T]$ it is at most $Ty/(\pi d(v)^2)$, where $d(v)\ge T/2$ is the distance from $v$ to $[T,2T]$. Since $G^+(v)\ll\log(2+|v|)$, these $v$ contribute $O(\log^2T)$. The remaining $v$ contribute at most $2\int_{T/2}^{3T}(\log^+|\zeta(\tfrac12+iv)|)^2\,dv\ll T\log\log T$, by the case $\sigma=\tfrac12$ applied to $[T/2,T]$, $[T,2T]$ and $[2T,4T]$.
\end{proof}

\begin{lemma}[Maximum to the right of the critical line]\label{lem:max}
For large $T$ put $L=\log T$ and $h=4\pi/L$, and for $t\in\mathbb R$ let
\[
m_T(t)=\sup\bigl\{\log|\zeta(\sigma+iv)|:\ \tfrac12\le\sigma\le3,\ |v-t|\le h\bigr\}.
\]
Then, uniformly for $A\ge1$,
\[
\operatorname{meas}\{t\in[T,2T]:\ m_T(t)>A\}\ll\frac{T\log\log T}{A^2}.
\]
\end{lemma}
\begin{proof}
We may assume that $A$ exceeds a large absolute constant. Put $a=\tfrac12-1/L$, $G(v)=\log|\zeta(a+iv)|$ and $\phi=G^+\mathbf 1_{[T/2,3T]}$. Let $t\in[T,2T]$, $\tfrac12\le\sigma\le3$ and $|v-t|\le h$, and put $y=\sigma-a\in[1/L,3]$. By \eqref{eq:poisson2},
\[
\log|\zeta(\sigma+iv)|\le(P_y*\phi)(v)+\int_{w\notin[T/2,3T]}P_y(v-w)G^+(w)\,dw+O(1/T).
\]
For $w\notin[T/2,3T]$ we have $|v-w|\ge T/3$, and $G^+(w)\ll\log(2+|w|)$, so the integral is $\ll y\log T/T\le1$. Since $|v-t|\le h=4\pi/L\le4\pi y$,
\[
y^2+(t-w)^2\le y^2+2(v-w)^2+2(v-t)^2\le c_1\bigl(y^2+(v-w)^2\bigr),\qquad c_1=1+32\pi^2 .
\]
Hence $P_y(v-w)\le c_1P_y(t-w)$, and $(P_y*\phi)(v)\le c_1(P_y*\phi)(t)\le c_1(\mathcal M\phi)(t)$, where $\mathcal M$ is the Hardy--Littlewood maximal operator. The last step holds because $P_y$ is even, decreasing in $|u|$ and has integral $1$ \cite[Ch.~III, \S2.2]{Stein}. Therefore $m_T(t)\le c_1(\mathcal M\phi)(t)+2$.

Put $b=(A-2)/c_1$ and $\psi=\phi\,\mathbf 1_{\{\phi>b/2\}}$. Then $\phi\le b/2+\psi$, so $\mathcal M\phi\le b/2+\mathcal M\psi$. By the weak type $(1,1)$ inequality for $\mathcal M$ \cite[Ch.~I, \S1]{Stein},
\[
\operatorname{meas}\{t\in[T,2T]:\ m_T(t)>A\}\le\operatorname{meas}\{\mathcal M\psi>b/2\}\ll\frac1b\int\psi\le\frac2{b^2}\int\phi^2 .
\]
By the functional equation and Lemma~\ref{lem:chi}, for $T/2\le w\le3T$, $G(w)=\log|\chi(a+iw)|+\log|\zeta(\tfrac12+\tfrac1L+iw)|$ and $\log|\chi(a+iw)|=L^{-1}\log(w/2\pi)+O(1/w)\le2$. So $G^+(w)\le2+\log^+|\zeta(\tfrac12+\tfrac1L+iw)|$ there, and Lemma~\ref{lem:moment} (with $\sigma=\tfrac12+\tfrac1L$, applied to $[T/2,T]$, $[T,2T]$ and $[2T,4T]$) gives $\int\phi^2\ll T\log\log T$. Since $b\asymp A$, the lemma follows.
\end{proof}

\section{Separators and the phase-slip counter}\label{sec:counter}

This section sets up the counter for a general target and identifies it with the phase of $\zeta$ at the balance abscissa. Let $f$ be analytic on $\{t>t_0\}$. We call real numbers $\sigma_L<\sigma_R$ with $\sigma_L<0$, together with $\varphi_0\in\mathbb R$, a \emph{frame} for $f$ if, as $t\to\infty$,
\begin{enumerate}
\item[(F1)] $|f/\zeta|(\sigma_L+it)\to0$ and $|\zeta/f|(\sigma_R+it)\to0$;
\item[(F2)] $\arg f(\sigma_R+it)\to\varphi_0$, for a continuous branch;
\item[(F3)] every solution of $\zeta(s)=f(s)$ with $t$ sufficiently large has $\sigma_L<\sigma<\sigma_R$.
\end{enumerate}
Given a frame, define the \emph{separators} $\tau_m$ by $\theta(\tau_m)=(m+\tfrac12)\pi-\tfrac12\varphi_0$, for large $m$; these are shifted Gram points. Equivalently, $\chi(\tfrac12+i\tau_m)=-e^{i\varphi_0}$: the separators are $a$-points of $\chi$ on the critical line, where $\zeta(\tfrac12+it)$ meets the line $e^{i(\varphi_0+\pi)/2}\mathbb R$ (cf.~\cite{KalpokasSteuding,SteudingSuriajaya}). Lester \cite{Lester} uses the name shifted Gram points for the solutions of $\theta(t)\equiv-\arg a\pmod\pi$; our normalisation is adapted instead to solutions off the line. Let $n_m$ be the number of solutions, counted with multiplicity, with $\tau_{m-1}<t\le\tau_m$, and let $N_f(T)$ be the number with $t_0<t\le T$. If a solution lies on a line $t=\tau_m$, we move $\tau_m$ up by at most $\tau_m^{-2}$; this changes $2\theta(\tau_m)$ by $O(\tau_m^{-2}\log\tau_m)$, which is absorbed in the $O(1/t)$ terms below. Put $\alpha_m=\arg\zeta(1-\sigma_L-i\tau_m)$, with the Dirichlet-series branch, which is available since $1-\sigma_L>1$.

We shall use four kinds of frames:
\begin{enumerate}
\item[(i)] $(\sigma_L,-\eta,\varphi_0)$ for the targets of Section~\ref{sec:big}, with $\varphi_0$ from (H2);
\item[(ii)] $(-1,2,\arg c-\kappa\pi/2)$ for $f(s)=c(-s)^{\kappa}$ with $c\ne0$ and $0<\kappa<\tfrac32$;
\item[(iii)] $(-1,2,0)$ for $f=\ell$;
\item[(iv)] $(-1,2,\varphi_0)$ for a target $f$ of Section~\ref{sec:slow} that satisfies \textup{(S3)}; case (iii) is the special case $f=\ell$, $\varphi_0=0$.
\end{enumerate}
In case (ii), (F1) and (F2) are immediate, since $|f/\zeta|(-1+it)\ll t^{\kappa-3/2}$ and $\arg(-2-it)\to-\pi/2$. Condition (F3) follows from Lemma~\ref{lem:left}, applied with $K=\kappa<\tfrac32$ and $\sigma_L=-1$, and from $|\zeta(\sigma+it)|\le\zeta(2)<|f(\sigma+it)|$ for $\sigma\ge2$. Cases (iii) and (iv) are similar. (F1) holds because $|f(-1+it)|=t^{o(1)}$ by \textup{(S1)} and $|f(2+it)|\to\infty$, (F2) is \textup{(S3)}, and (F3) follows from the first assertion of Lemma~\ref{lem:strip}, which needs no restriction on $K$. Case (i) is checked in the proof of Theorem~\ref{thm:gram}.

\begin{proposition}[Phase-slip counter]\label{prop:counter}
Let $(\sigma_L,\sigma_R,\varphi_0)$ be a frame for $f$. There are $m_0$, a constant $C_f$ and integers $k_m$ $(m\ge m_0)$ such that the change of $\arg(\zeta-f)$ along the separator from $\sigma_L+i\tau_m$ to $\sigma_R+i\tau_m$ equals $-\alpha_m+2\pi k_m+o(1)$, and
\[
n_m=1+k_{m-1}-k_m,\qquad N_f(\tau_m)=m+C_f-k_m .
\]
\end{proposition}
\begin{proof}
Apply the argument principle to the rectangle $[\sigma_L,\sigma_R]\times[\tau_{m-1},\tau_m]$. By (F1) and (F3), no solution lies on its vertical sides when $m$ is large.
\begin{itemize}
\item Left side. Here $\zeta-f=\chi(s)\zeta(1-s)(1-f/\zeta)$ with $|f/\zeta|\to0$, and $\arg\zeta(1-s)$ is given by the Dirichlet series. Going down from $\tau_m$ to $\tau_{m-1}$, Lemma~\ref{lem:chi} shows that $\arg\chi$ changes by $2\theta(\tau_m)-2\theta(\tau_{m-1})+o(1)=2\pi+o(1)$. So $\arg(\zeta-f)$ changes by $2\pi+\alpha_{m-1}-\alpha_m+o(1)$.
\item Right side. Here $\zeta-f=-f(1-\zeta/f)$ with $|\zeta/f|\to0$, so by (F2) the change is $o(1)$.
\item Separators. By Lemma~\ref{lem:chi}, $\arg\chi(\sigma+i\tau_m)\equiv-2\theta(\tau_m)\equiv\varphi_0+\pi\pmod{2\pi}$ up to $O(1/\tau_m)$, uniformly for $\sigma\in[\sigma_L,\sigma_R]$. Hence $\arg(\zeta-f)\equiv\pi+\varphi_0+\alpha_m+o(1)$ at $\sigma_L$ and $\arg(\zeta-f)\equiv\pi+\varphi_0+o(1)$ at $\sigma_R$. So the change along the separator is $-\alpha_m+2\pi k_m+o(1)$ for an integer $k_m$.
\end{itemize}
Traversing the rectangle positively, the total change is $2\pi+2\pi(k_{m-1}-k_m)+o(1)$, and it equals $2\pi n_m$. This gives the first identity for large $m$, and telescoping gives the second.
\end{proof}

\begin{corollary}[Block criterion]\label{cor:block}
Gram's law ($n_m=1$) holds for all $m$ in a block $m_1<m\le m_2$ if and only if $k_{m_1}=k_{m_1+1}=\dots=k_{m_2}$. Moreover $N_f(\tau_m)-(\theta(\tau_m)+\varphi_0/2)/\pi=C_f-\tfrac12-k_m$. Hence, if Gram's law holds for all large $m$, then $N_f(T)-(\theta(T)+\varphi_0/2)/\pi$ is bounded; and if $k_m$ is unbounded, then Gram's law fails for infinitely many $m$.
\end{corollary}

The counter is determined by the phase of $\zeta$ at the balance abscissa, whenever the modulus ratio crosses $1$ only once along the separator.

\begin{proposition}[Identification of the counter]\label{prop:ident}
Let $(\sigma_L,\sigma_R,\varphi_0)$ be a frame for $f$, let $m$ be large, and write $\tau=\tau_m$ and $R(\sigma)=|\zeta(\sigma+i\tau)/f(\sigma+i\tau)|$. Suppose there is $\sigma_c\in(\sigma_L,\sigma_R)$ with $\sigma_c<\tfrac12$, $R\ge1$ on $[\sigma_L,\sigma_c]$ and $R\le1$ on $[\sigma_c,\sigma_R]$. Put $x_c=1-\sigma_c$, and let $\varepsilon(\sigma)$ be the continuous branch of $\arg f(\sigma+i\tau)-\varphi_0$ on $[\sigma_L,\sigma_R]$ with $\varepsilon(\sigma_R)=o(1)$. Then $\zeta(x+i\tau)\ne0$ for $x\ge x_c$, and $k_m$ is the integer nearest to $\psi_m/2\pi$, where
\[
\psi_m=-\pi S(\tau;x_c)-\varepsilon(\sigma_c)+r_m,\qquad |r_m|\le C/\tau ,
\]
and $C$ depends only on $\sigma_L$ and $\sigma_R$. Moreover $\psi_m\notin\pi+2\pi\mathbb Z$.
\end{proposition}
\begin{proof}
On the separator, Lemma~\ref{lem:chi} gives $\Im\log\chi(\sigma+i\tau)=-2\theta(\tau)+r(\sigma)$ with $|r|\le C/\tau$, and $-2\theta(\tau)=-(2m+1)\pi+\varphi_0$. So $\chi(s)=-|\chi(s)|e^{i(\varphi_0+r)}$ and
\[
\zeta-f=-f\,(1+u),\qquad u=-\frac{\zeta}{f}=|\chi(s)|\,e^{i(\varphi_0+r)}\,\frac{\zeta(1-s)}{f(s)} .
\]
On $[\sigma_L,\sigma_c]$ we have $R\ge1$, so $\zeta(s)\ne0$ and hence $\zeta(1-s)\ne0$; for $x>1$ there are no zeros anyway. Hence a continuous argument of $u$ on $[\sigma_L,\sigma_c]$ is
\[
\psi(\sigma)=r(\sigma)+\Im\log\zeta(1-\sigma-i\tau)-\varepsilon(\sigma),
\]
where $\log\zeta(1-\sigma-i\tau)$ is continued from $\sigma=\sigma_L$, at which point it is given by the Dirichlet series. By conjugation symmetry, $\Im\log\zeta(x-i\tau)=-\pi S(\tau;x)$ when both are continued from $\Re s>1$. So $\psi(\sigma_c)=\psi_m$ with $r_m=r(\sigma_c)$, and $\psi(\sigma_L)=\alpha_m+r(\sigma_L)-\varepsilon(\sigma_L)$.

Now apply Lemma~\ref{lem:cross} on $[\sigma_L,\sigma_R]$ with $\sigma^*=\sigma_c$. Its hypotheses hold, because $\zeta-f$ does not vanish on the separator. By (F1), $|u(\sigma_R)|\to0$ and $|u(\sigma_L)|^{-1}\to0$. So the change of $\arg(1+u)$ along the separator is $-\psi(\sigma_L)+2\pi j+o(1)$, where $j$ is the integer nearest to $\psi_m/2\pi$, and $\psi_m\notin\pi+2\pi\mathbb Z$. The change of $\arg(-f)$ is $\varepsilon(\sigma_R)-\varepsilon(\sigma_L)$. Adding, the change of $\arg(\zeta-f)$ is $-\alpha_m+2\pi j+o(1)$. Comparing with Proposition~\ref{prop:counter} gives $k_m=j$.
\end{proof}

The single-crossing hypothesis is verified unconditionally when $x_c>1$ (Step~2 of the proof of Theorem~\ref{thm:sharp}), and at almost all separators for the targets of Sections~\ref{sec:smallkappa} and~\ref{sec:gramlog} (Lemma~\ref{lem:ident}, with Step~4 of the proof of Theorem~\ref{thm:smallkappa} and Step~3 of the proof of Theorem~\ref{thm:gramfail}). In each case $\varepsilon=o(1)$, so the counter is the integer nearest to $-\tfrac12S(\tau_m;x_c)+o(1)$, as announced in the introduction.

\section{Targets growing faster than \texorpdfstring{$t^{1/2}$}{t\^{}(1/2)}: the regime \texorpdfstring{$x_c>1$}{xc>1}}\label{sec:big}

In this section the balance abscissa satisfies $x_c>1$, where the phase of $\zeta$ is controlled by the absolutely convergent Euler product, and everything is unconditional.

\begin{theorem}\label{thm:gram}
Let $\kappa>\kappa_1$ and $\varphi_0\in\mathbb R$. Choose $\eta>0$ with $\tfrac12+\eta<\kappa$ and $A(1+\eta)<\pi$ (possible since $\kappa>\kappa_1$). Let $f$ be analytic on $\{t>t_0\}$ with $|f(s)|\le C_1|s|^{K}$ there, and fix $\sigma_L<\min(-\eta,\tfrac12-K)$. Suppose that, as $t\to\infty$,
\begin{enumerate}
\item[(H1)] $\log|f(\sigma+it)|=\kappa\log t+o(\log t)$ uniformly for $\sigma\in[\sigma_L,-\eta]$, and $|f(\sigma+it)|\ge t^{\kappa-o(1)}$ uniformly for $\sigma\ge-\eta$;
\item[(H2)] $\arg f(\sigma+it)=\varphi_0+o(1)$ uniformly for $\sigma\in[\sigma_L,-\eta]$, for a continuous branch.
\end{enumerate}
Define $\tau_m$ by $\theta(\tau_m)=(m+\tfrac12)\pi-\tfrac12\varphi_0$. Then there is $m_0$ such that, for every $m\ge m_0$, the strip $\tau_{m-1}<t\le\tau_m$ contains exactly one solution $\rho_m$ of $\zeta(s)=f(s)$. Moreover $\rho_m$ is simple, $\sigma_L<\beta_m<-\eta$, and
\[
\Bigl|\theta(\gamma_m)-m\pi+\tfrac12\varphi_0\Bigr|\le\tfrac12A(1+\eta)+o(1)<\tfrac{\pi}{2}.
\]
\end{theorem}

\begin{proof}
Write $L=A(1+\eta)<\pi$. All $o(1)$ terms are as $m\to\infty$, uniformly on the segments in question.

\emph{Where solutions can be.} For $\sigma\le\sigma_L$, Lemma~\ref{lem:left} gives $|\zeta|>|f|$. For $\sigma\ge-\eta$, we have $|\zeta(\sigma+it)|\ll t^{1/2+\eta}\log t$: use the functional equation on $[-\eta,0]$, convexity on $[0,1]$, and $|\zeta|\ll\log t$ for $\sigma\ge1$. By (H1) and $\kappa>\tfrac12+\eta$, $|\zeta|<|f|$ there for large $t$. So every solution with large $t$ has $\sigma_L<\beta<-\eta$. Together with (H2), this shows that $(\sigma_L,-\eta,\varphi_0)$ is a frame for $f$, so Proposition~\ref{prop:counter} applies.

\emph{Phase of $\zeta$ for $\sigma\in[\sigma_L,-\eta]$.} Here $\Re(1-s)\ge1+\eta$, so $|\arg\zeta(1-s)|\le L$ by Lemma~\ref{lem:arcsin}(a). By Lemma~\ref{lem:chi},
\[
\zeta(s)=|\zeta(s)|\exp\bigl(i(-2\theta(t)+\alpha(s))\bigr),\qquad \alpha(s)=\arg\zeta(1-s)+O(1/t),\quad|\alpha|\le L+o(1).
\]

\emph{The counter vanishes.} On $t=\tau_m$ we have $-2\theta(\tau_m)\equiv\varphi_0+\pi\pmod{2\pi}$, hence
\[
\zeta(s)-f(s)=-e^{i\varphi_0}\Bigl(|\zeta(s)|e^{i\alpha(s)}+|f(s)|e^{i\varepsilon(s)}\Bigr),\qquad\varepsilon=o(1)\ \text{by (H2)}.
\]
The two vectors in the bracket make an angle of at most $L+o(1)<\pi$. So the bracket never vanishes, and its argument stays in the interval between $\alpha$ and $\varepsilon$, which has length less than $\pi$. At $\sigma=\sigma_L$ the argument is $\alpha+o(1)$, since $|f/\zeta|\to0$, and at $\sigma=-\eta$ it is $o(1)$, since $|\zeta/f|\to0$. Hence the change of $\arg(\zeta-f)$ along the separator, from left to right, is $-\alpha(\sigma_L+i\tau_m)+o(1)=-\alpha_m+o(1)$. So $k_m=0$ for all large $m$, and $n_m=1$ by Proposition~\ref{prop:counter}. Since $n_m$ counts solutions with multiplicity, the solution $\rho_m$ is simple.

\emph{Localisation.} By \eqref{eq:phase}, $2\theta(\gamma_m)\equiv\alpha(\rho_m)-\varphi_0+o(1)\pmod{2\pi}$, with $|\alpha|\le L+o(1)$. Since $\theta(\tau_{m-1})<\theta(\gamma_m)\le\theta(\tau_m)$, this forces $\theta(\gamma_m)-m\pi+\tfrac12\varphi_0=\tfrac12\alpha(\rho_m)+o(1)$.
\end{proof}

\begin{corollary}[Counting]\label{cor:count}
Under the hypotheses of Theorem~\ref{thm:gram}, let $N_f(T)$ count the solutions with $t_0<t\le T$. Then there is a constant $C_f$ with
\[
N_f(T)=\frac{\theta(T)+\varphi_0/2}{\pi}+C_f+E(T),\qquad |E(T)|\le1\quad(T\ge T_0).
\]
Hence $N_f(T)-N(T)=-S(T)+O(1)$. Since $S(T)$ is unbounded \cite[Ch.~IX]{Titchmarsh}, so is $N_f(T)-N(T)$.
\end{corollary}
\begin{proof}
By the proof of Theorem~\ref{thm:gram}, $k_m=0$ for large $m$, so $N_f(\tau_m)=m+C$ by Proposition~\ref{prop:counter}. For $\tau_{m-1}<T\le\tau_m$ we have $N_f(T)\in\{N_f(\tau_{m-1}),N_f(\tau_m)\}$ and $(\theta(T)+\varphi_0/2)/\pi\in(m-\tfrac12,m+\tfrac12]$; take $C_f=C-\tfrac12$.
\end{proof}

\begin{corollary}[Polynomial targets]\label{cor:poly}
Let $P$ be a polynomial of degree $d\ge1$ with leading coefficient $c$. Then, for all large $m$, exactly one solution of $\zeta(s)=P(s)$ lies in each strip $\tau_{m-1}<t\le\tau_m$, where $\varphi_0=\arg c+d\pi/2$. Each such solution is simple, and
\[
\beta=\tfrac12-\frac{d\log\gamma}{\log(\gamma/2\pi)}+O\Bigl(\frac1{\log\gamma}\Bigr).
\]
For the fixed points ($P(s)=s$) this gives $\beta\to-\tfrac12$, and the number of non-real fixed points with $0<t\le T$ is $\frac{T}{2\pi}\log\frac{T}{2\pi e}+O(1)$. This is consistent with, and localises, the disc count of \cite{LiSteuding}.
\end{corollary}
\begin{proof}
Take $\kappa=K=d$ and $\eta=\tfrac12\bigl((\sigma_1-1)+(d-\tfrac12)\bigr)$ in Theorem~\ref{thm:gram}, so that $\sigma_1-1<\eta<d-\tfrac12$. The formula for $\beta$ follows from \eqref{eq:modulus}, because $|\zeta(1-\rho)|$ is bounded above and below. Only finitely many non-real fixed points have $0<t\le t_0$: for $\sigma\le-M$ they can only lie in small discs about the trivial zeros, where Rouch\'e's theorem gives one solution, real by symmetry.
\end{proof}

\subsection{Sharpness of $\kappa_1$}

\begin{theorem}[Sharpness]\label{thm:sharp}
Let $\tfrac12<\kappa<\kappa_1$ and $\varphi_0\in\mathbb R$. Let $f$ be analytic on $\{t>t_0\}$ with $|f(s)|\le C_1|s|^K$, fix $\eta\in(0,\kappa-\tfrac12)$, and fix $\sigma_L<\min(-\eta,\ \tfrac12-K)$. Assume (H2), and
\begin{enumerate}
\item[(H1$'$)] $\log|f(\sigma+it)|=\kappa\log t+O(1)$ and $\partial_\sigma\log|f(\sigma+it)|=o(\log t)$, uniformly for $\sigma\in[\sigma_L,-\eta]$, and $|f(\sigma+it)|\ge t^{\kappa-o(1)}$ uniformly for $\sigma\ge-\eta$.
\end{enumerate}
With $\tau_m$ as in Theorem~\ref{thm:gram}, let $n_m$ be the number of solutions of $\zeta(s)=f(s)$ with $\tau_{m-1}<t\le\tau_m$. Then
\[
\#\{m:\ \tau_m\le T,\ n_m\neq1\}\gg T.
\]
\end{theorem}

\begin{corollary}[Dichotomy]\label{cor:dich}
Let $f(s)=c(-s)^{\kappa}$ with $c\neq0$, $\kappa>\tfrac12$ and $\kappa\neq\kappa_1$. Then $n_m=1$ for all large $m$ if and only if $\kappa>\kappa_1$.
\end{corollary}
\begin{proof}
Here $\log|f|=\log|c|+\kappa\log|s|$, $\partial_\sigma\log|f|=\kappa\sigma/|s|^2$, and $\arg f\to\arg c-\kappa\pi/2$. Apply Theorems~\ref{thm:gram} and~\ref{thm:sharp}.
\end{proof}

\begin{proof}[Proof of Theorem~\ref{thm:sharp}]
\emph{Step 1: the counter.} The first step of the proof of Theorem~\ref{thm:gram} uses only $\kappa>\tfrac12+\eta$, $\sigma_L<\tfrac12-K$ and the lower bound for $|f|$ on $\sigma\ge-\eta$, which \textup{(H1$'$)} provides, not any bound on $\arg\zeta$. So $(\sigma_L,-\eta,\varphi_0)$ is again a frame for $f$, and Proposition~\ref{prop:counter} gives integers $k_m$ with $n_m=1+k_{m-1}-k_m$. Thus $n_m\neq1$ if and only if $k_{m-1}\neq k_m$.

\emph{Step 2: identifying $k_m$.} Put $R=|\zeta/f|$ on the separator. By Lemma~\ref{lem:chi} and (H1$'$),
\[
\log R=(\tfrac12-\sigma-\kappa)\log t+O(1),\qquad \partial_\sigma\log R=-(1+o(1))\log t,
\]
since $\partial_\sigma\log|\zeta(1-s)|$ is bounded for $\Re(1-s)\ge1+\eta$. Hence $R$ is strictly decreasing in $\sigma$ and equals $1$ at a unique $\sigma_c(\tau_m)=\tfrac12-\kappa+O(1/\log\tau_m)$. So Proposition~\ref{prop:ident} applies. Since $x_c=1-\sigma_c>1+\eta$, we have $\pi S(\tau_m;x_c)=\Im\log\zeta(x_c+i\tau_m)$ with the Dirichlet-series branch, and $\varepsilon=o(1)$ by (H2). Therefore, for all large $m$,
\begin{equation}\label{eq:kround}
k_m=\Bigl\lfloor\frac{\psi_m+\pi}{2\pi}\Bigr\rfloor,\qquad\psi_m=-\Im\log\zeta(x_c+i\tau_m)+o(1);
\end{equation}
that is, $k_m$ is the integer nearest to $\psi_m/2\pi$, where $\psi_m\notin\pi+2\pi\mathbb Z$ is the number of Proposition~\ref{prop:ident}.

\emph{Step 3: phase slips.} Since $\sigma_c=\tfrac12-\kappa+O(1/\log\tau)$,
\[
\psi_m=-\Im\log\zeta(\tfrac12+\kappa+i\tau_m)+o(1).
\]
Since $\tfrac12+\kappa<\sigma_1$, we have $A(\tfrac12+\kappa)>\pi$. Fix $\delta\in(0,\tfrac{\pi}{4})$ with $\pi+2\delta<A(\tfrac12+\kappa)$. Apply Lemma~\ref{lem:arcsin}(b) at $\sigma_0=\tfrac12+\kappa$ with $v=-(\pi+2\delta)$ and with $v=0$ (tolerance $\delta/2$). This yields two families of intervals, $\mathcal I$ and $\mathcal J$, each of length $d$, such that every interval of length $L$ contains a member of each family. For large $m$:
\begin{itemize}
\item if $\tau_m$ lies in a member of $\mathcal I$, then $\psi_m\in(\pi+\delta,\pi+3\delta)$, so $k_m=1$;
\item if $\tau_m$ lies in a member of $\mathcal J$, then $|\psi_m|<\delta$, so $k_m=0$.
\end{itemize}
The spacing of the $\tau_m$ tends to $0$, so every member of $\mathcal I\cup\mathcal J$ contains some $\tau_m$. Between a member of $\mathcal I$ and the next member of $\mathcal J$, some consecutive pair has $k_{m-1}\neq k_m$, that is, $n_m\neq1$. There are $\gg T/L$ such alternations in $[T,2T]$.
\end{proof}

\begin{remark}
The proof shows that failures occur in phase slips. Gram's law fails exactly when the counter $k_m=\lfloor(\psi_m+\pi)/2\pi\rfloor$ jumps, that is, essentially when $\arg\zeta(\tfrac12+\kappa+it)$ crosses an odd multiple of $\pi$ (see Step~5 of the proof of Theorem~\ref{thm:order}). Such crossings occur $\asymp T$ times in $[T,2T]$, against $\asymp T\log T$ Gram intervals, so in this regime one expects a failure proportion of order $1/\log T$; Theorem~\ref{thm:order} below shows that this is correct. For the zeta zeros, by contrast, the proportion is positive \cite{Trudgian}. Because $\arg\zeta(\tfrac12+\kappa+it)$ must exceed $\pi$, which needs the factors $p^{-it}$ for many primes $p$ to be aligned at once, the first slip for $\kappa\in(\tfrac12,\kappa_1)$ lies far beyond computational range.
\end{remark}

The upper bound in the following proposition is weaker than the heuristic $O(T)$, but it is enough to show that failures are rare.

\begin{proposition}[Failures have density zero]\label{prop:upper}
Under the hypotheses of Theorem~\ref{thm:sharp},
\[
\#\{m:\ \tau_m\le T,\ n_m\neq1\}=o(T\log T).
\]
Thus the proportion of failing intervals tends to $0$. By Theorem~\ref{thm:sharp} it is $\gg1/\log T$.
\end{proposition}
\begin{proof}
Put $x_0=\tfrac12+\kappa\in(1,\sigma_1)$ and $\Phi(t)=\Im\log\zeta(x_0+it)$, with the Dirichlet-series branch. We use the notation of the proof of Theorem~\ref{thm:sharp}.

\emph{Step 1: $\psi_m=-\Phi(\tau_m)+o(1)$.} By \eqref{eq:kround}, $\psi_m=-\Im\log\zeta(x_c+i\tau_m)+o(1)$, where $x_c=1-\sigma_c=x_0+O(1/\log\tau_m)$. For $x\ge1+\eta$ we have $|\partial_x\Im\log\zeta(x+i\tau)|\le|\zeta'/\zeta(1+\eta)|$.

\emph{Step 2: a failure forces $\Phi$ to be close to $\pi+2\pi\mathbb Z$.} If $n_m\ne1$, then $k_{m-1}\ne k_m$. By \eqref{eq:kround}, the closed interval between $\psi_{m-1}$ and $\psi_m$ then contains a point of $\pi+2\pi\mathbb Z$. Now $|\Phi'(t)|\le|\zeta'/\zeta(x_0)|$ and $\tau_m-\tau_{m-1}\to0$, so Step~1 gives $|\psi_m-\psi_{m-1}|=o(1)$. Hence for every $\eta_0>0$ there is $m(\eta_0)$ such that, if $m\ge m(\eta_0)$ and $n_m\ne1$, then $\operatorname{dist}(\Phi(t),\pi+2\pi\mathbb Z)\le\eta_0$ for all $t\in[\tau_{m-1},\tau_m]$.

\emph{Step 3: a limit law.} Put $Y=Y_{x_0}$ (Section~\ref{sec:main}). The series converges absolutely, and $|Y|\le A(x_0)$. Since $x_0>1$, the series for $\Phi$ converges uniformly. Truncating at $p\le y$ and applying the Kronecker--Weyl theorem to $(t\log p)_{p\le y}$, whose frequencies are linearly independent over $\mathbb Q$, gives
\[
\frac1T\operatorname{meas}\{t\in[T,2T]:\Phi(t)\in U\}\longrightarrow\mathbb P(Y\in U)
\]
for every finite union $U$ of intervals whose endpoints are not atoms of $Y$ (this is the case $x_0>1$ of the Bohr--Jessen theorem \cite{BohrJessen}). In fact $Y$ has no atoms. Write $Y=Y_2+Y'$, where $Y_2=-\arg(1-2^{-x_0}e^{i\vartheta_2})$ is independent of $Y'$. Then $Y_2$ is a non-constant real-analytic function of $\vartheta_2$, so each of its level sets is finite, and $Y_2$ has no atoms; hence neither does $Y$. Put $U_{\eta_0}=\{y\in[-A(x_0),A(x_0)]:\operatorname{dist}(y,\pi+2\pi\mathbb Z)\le\eta_0\}$. Then $\mathbb P(Y\in U_{\eta_0})\to0$ as $\eta_0\to0$.

\emph{Step 4: counting.} The intervals $[\tau_{m-1},\tau_m]$ with $T<\tau_m\le2T$ are disjoint. Since $\theta'(t)=\tfrac12\log(t/2\pi)+O(t^{-2})$, each has length at least $(2\pi-o(1))/\log T$. By Step~2, for $m\ge m(\eta_0)$ the failing ones lie in $\{t\in[T-1,2T]:\Phi(t)\in U_{\eta_0}\}$. By Step~3,
\[
\#\{m:\ T<\tau_m\le2T,\ n_m\ne1\}\le\frac{\log T}{2\pi}\bigl(\mathbb P(Y\in U_{\eta_0})+o(1)\bigr)T+O_{\eta_0}(1).
\]
Let $T\to\infty$ and then $\eta_0\to0$, and sum over dyadic ranges.
\end{proof}

For power targets the upper bound can be brought down to the order of the lower bound. The proof counts phase slips directly, by Jensen's formula, so that no information about the distribution of $\arg\zeta$ at scale $1/\log T$ is needed.

\begin{theorem}[Order of the failure count]\label{thm:order}
Let $f(s)=c(-s)^{\kappa}$ with $c\ne0$ and $\tfrac12<\kappa<\kappa_1$. Let $x_0=\tfrac12+\kappa$ and $Y=Y_{x_0}$.
\begin{enumerate}
\item[(a)] $\displaystyle\sum_{m:\,\tau_m\le T}|n_m-1|\asymp T$. In particular $\#\{m:\tau_m\le T,\ n_m\ne1\}\asymp T$, so the proportion of failing intervals is $\asymp1/\log T$.
\item[(b)] For every $j\in\mathbb Z$,
\[
\lim_{M\to\infty}\frac1M\#\{m\le M:\ k_m=j\}=\mathbb P\bigl((2j-1)\pi<-Y<(2j+1)\pi\bigr).
\]
In particular $k_m\ne0$ for a proportion $\mathbb P(|Y|>\pi)>0$ of all $m$.
\end{enumerate}
\end{theorem}

Since $N_f(\tau_m)=m+C-k_m$ by telescoping, part (b) says that for $\tfrac12<\kappa<\kappa_1$ the solutions stay attached to the shifted Gram points, up to a random offset of a whole number of intervals. The offset changes only $\asymp T$ times up to height $T$, yet it is non-zero at a positive proportion of the separators.

\begin{proof}
The lower bound in (a) is Theorem~\ref{thm:sharp}, since $\sum_m|n_m-1|\ge\#\{m:n_m\ne1\}$. We use the notation of the proofs of Theorem~\ref{thm:sharp} and Proposition~\ref{prop:upper}. Here $\varepsilon(s)=\arg f(s)-\varphi_0=-\kappa\arctan(\sigma/t)$.

\emph{Step 1: an interpolating phase.} Step~2 of the proof of Theorem~\ref{thm:sharp} applies at every large height $t$, not only at separators. So $\sigma\mapsto\log|\zeta/f|(\sigma+it)$ has a unique zero $\sigma_c(t)\in(\sigma_L,-\eta)$, which is a real-analytic function of $t$ with $\sigma_c(t)=1-x_0+O(1/\log t)$. Put
\[
\Upsilon(t)=\alpha\bigl(\sigma_c(t)+it\bigr)-\varepsilon\bigl(\sigma_c(t)+it\bigr),\qquad \alpha(\sigma+it)=\Im\log\zeta(1-\sigma-it)+\Im\log\chi(\sigma+it)+2\theta(t),
\]
with the branch of $\log\chi$ from Lemma~\ref{lem:chi} and the Dirichlet-series branch of $\log\zeta$. Then $\Upsilon(\tau_m)=\psi_m$, the number of Proposition~\ref{prop:ident}, so by \eqref{eq:kround}
\begin{equation}\label{eq:khalf}
\bigl|k_m-\Upsilon(\tau_m)/2\pi\bigr|\le\tfrac12\qquad(m\ \text{large}).
\end{equation}
(If $\tau_m$ was moved as in Section~\ref{sec:counter}, let the displacement tend to $0$; $k_m$ does not change, and $\psi_m\to\Upsilon(\tau_m)$.) As in Step~1 of the proof of Proposition~\ref{prop:upper}, $\Upsilon(t)=-\Phi(t)+O(1/\log t)$, where $\Phi(t)=\Im\log\zeta(x_0+it)$.

\emph{Step 2: analytic continuation.} Let $\nu=(x_0-1)/4$ and $S=\{z:\ |\Im z|<\nu,\ \Re z>t_2\}$, with $t_2$ large. We claim that $\Upsilon$ extends to an analytic function $\tilde\Upsilon$ on $S$ with $|\tilde\Upsilon|\le M:=\log\zeta(1+\nu)+1$.

Let $F(s)=\log\chi(s)+\log\zeta(1-s)-\log f(s)$ for $|\Re s-(1-x_0)|<3\nu$ and $\Im s>t_2-1$. There $\Re(1-s)>1+\nu$. Put $F^*(w)=\overline{F(\bar w)}$. For $|\sigma-(1-x_0)|<2\nu$, $|\Im\sigma|<1$ and $z\in S$, define
\[
G(\sigma,z)=\tfrac12\bigl[F(\sigma+iz)+F^*(\sigma-iz)\bigr],\qquad P(\sigma,z)=\tfrac1{2i}\bigl[F(\sigma+iz)-F^*(\sigma-iz)\bigr]+2\theta(z)+\varphi_0,
\]
where $\theta(z)$ is the analytic continuation of $\theta$. For real $\sigma$ and $z$ we have $G=\log|\zeta/f|(\sigma+iz)$ and $P=\alpha-\varepsilon$ at $\sigma+iz$. Hence $\Upsilon(t)=P(\sigma_c(t),t)$.

By Stirling's formula, $\chi'/\chi(w)=-\log(w/2\pi i)+O(1/|w|)$ uniformly for $|\Re w|\le2$, $\Im w\ge1$. Also $\zeta'/\zeta(1-s)$ and its derivatives are bounded on the region above, and $f'/f(s)=\kappa/s$. From this, uniformly in the domain of $(\sigma,z)$:
\begin{enumerate}
\item[(i)] $\partial_\sigma G=-\log(z/2\pi)+O(1)$, $\partial_zG=O(1)$, and all second derivatives of $G$ are $O(1)$;
\item[(ii)] $|P(\sigma,z)|\le\log\zeta(1+\nu)+O(1/\Re z)$.
\end{enumerate}
For (ii), note that $2\theta(z)=\frac1{2i}\bigl[\log\chi^*(\tfrac12-iz)-\log\chi(\tfrac12+iz)\bigr]$: this holds for real $z$ by the choice of branch, hence for all $z$. So the $\chi$-terms of $P$ equal $\frac1{2i}\int_{1/2}^{\sigma}\bigl[\chi'/\chi(u+iz)-\chi^{*\prime}/\chi^*(u-iz)\bigr]du=O(1/|z|)$. The $\zeta$-terms are $\frac1{2i}[\log\zeta(1-\sigma-iz)-\log\zeta(1-\sigma+iz)]$, and both arguments have real part at least $1+\nu$. The $f$-terms equal $\kappa\arctan(\sigma/z)=O(1/|z|)$.

Now fix a large real $t$ and let $z=t+iy$ with $|y|<\nu$. By (i), for $|\sigma-\sigma_c(t)|=\rho\le\nu$ we have $|G(\sigma,t)|\ge\rho(\tfrac12\log t-C\rho)\ge\tfrac14\rho\log t$ and $|G(\sigma,z)-G(\sigma,t)|\le C|y|$. By Rouch\'e's theorem, $G(\cdot,z)$ has exactly one zero $\sigma_c(z)$ in $|\sigma-\sigma_c(t)|<\nu$, and $|\sigma_c(z)-\sigma_c(t)|\le4C\nu/\log t$. The formula
\[
\sigma_c(z)=\frac1{2\pi i}\oint_{|\sigma-\sigma_c(t)|=\nu}\sigma\,\frac{\partial_\sigma G(\sigma,z)}{G(\sigma,z)}\,d\sigma
\]
shows that $\sigma_c$ is analytic. So $\tilde\Upsilon(z)=P(\sigma_c(z),z)$ is analytic on $S$, it agrees with $\Upsilon$ on the real axis, and $|\tilde\Upsilon|\le M$ by (ii).

\emph{Step 3: $\Phi$ is nowhere flat.} Let $d=\nu/4$. We claim there is $\delta>0$ such that $\max_J\Phi-\min_J\Phi\ge\delta$ for every interval $J$ of length $d$. Let $\mathbb T=\prod_p\mathbb R/2\pi\mathbb Z$, which is compact, and for $\vartheta\in\mathbb T$ put
\[
\Phi_\vartheta(t)=-\sum_p\arg\bigl(1-p^{-x_0}e^{-i(t\log p+\vartheta_p)}\bigr).
\]
Since $x_0>1$ and $|\arg(1-w)|\le\tfrac\pi2|w|$ for $|w|<1$, the series converges uniformly in $(\vartheta,t)$, and also for complex $t$ with $|\Im t|\le\nu$, where the terms are $O(p^{\nu-x_0})$ and $x_0-\nu>1$. So $\Phi_\vartheta$ is real-analytic and $\vartheta\mapsto\Phi_\vartheta|_{[0,d]}$ is continuous into $C[0,d]$. Moreover $\Phi(t_0+t)=\Phi_\vartheta(t)$ with $\vartheta_p=t_0\log p$. Hence it suffices that the oscillation of $\Phi_\vartheta$ on $[0,d]$, a continuous function of $\vartheta$, never vanishes. If it did vanish, $\Phi_\vartheta$ would be constant on $[0,d]$, hence on $\mathbb R$ by analyticity. But expanding each term in powers of $p^{-x_0}$ and using uniform convergence,
\[
\lim_{U\to\infty}\frac1U\int_0^U\Phi_\vartheta(t)e^{it\log2}\,dt=\frac{2^{-x_0}e^{-i\vartheta_2}}{2i}\ne0,
\]
which is impossible for a constant function.

\emph{Step 4: counting phase slips.} Let $J\subset[T-1,2T]$ be an interval of length $d$, and let $c\in\pi+2\pi\mathbb Z$ with $|c|\le M$. By Steps~1 and~3, $\Upsilon$ oscillates by at least $\delta-O(1/\log T)\ge\delta/2$ on $J$, so there is $t^*\in J$ with $|\Upsilon(t^*)-c|\ge\delta/4$. Apply Jensen's formula to $\tilde\Upsilon-c$ on the disc $|z-t^*|\le2d$, which lies in $S$. There $|\tilde\Upsilon-c|\le2M$. So the number of zeros of $\tilde\Upsilon-c$ in $|z-t^*|\le d$, which contains $J$, is at most $\log(8M/\delta)/\log2$. Since $|\Upsilon|\le M$ on the real axis, only $O(M)$ values of $c$ occur. Covering $[T-1,2T]$ by $O(T)$ such intervals gives
\[
\#\{t\in[T-1,2T]:\ \Upsilon(t)\in\pi+2\pi\mathbb Z\}\ll T .
\]

\emph{Step 5: conclusion of (a).} Suppose $k_m>k_{m-1}$. By \eqref{eq:khalf}, $\Upsilon(\tau_{m-1})\le2\pi k_{m-1}+\pi$ and $\Upsilon(\tau_m)\ge2\pi k_m-\pi$. So each of the $k_m-k_{m-1}$ values $(2k+1)\pi$, $k_{m-1}\le k<k_m$, is taken by $\Upsilon$ on $[\tau_{m-1},\tau_m]$. The case $k_m<k_{m-1}$ is the same. Since $|n_m-1|=|k_{m-1}-k_m|$, and each point of $[T-1,2T]$ lies in at most two of the intervals $[\tau_{m-1},\tau_m]$,
\[
\sum_{T<\tau_m\le2T}|n_m-1|\le2\,\#\{t\in[T-1,2T]:\ \Upsilon(t)\in\pi+2\pi\mathbb Z\}\ll T .
\]
Summing over dyadic ranges proves (a).

\emph{Proof of (b).} By \eqref{eq:kround} and Step~1, $k_m=j$ if $-\Phi(\tau_m)\in[(2j-1)\pi+\eta_0,(2j+1)\pi-\eta_0]$, and $k_m=j$ implies $-\Phi(\tau_m)\in((2j-1)\pi-\eta_0,(2j+1)\pi+\eta_0)$, for any $\eta_0>0$ and $m\ge m(\eta_0)$. On a Gram interval $\Phi$ varies by $O(1/\log T)$, and the Gram intervals near height $T$ have length $(2\pi+o(1))/\log T$. So the number of $m$ with $T<\tau_m\le2T$ and $-\Phi(\tau_m)$ in an interval $V$ lies between $\frac{\log T}{2\pi}(1+o(1))\operatorname{meas}\{t\in[T,2T]:-\Phi(t)\in V^{\mp}\}+O(1)$, where $V^{\mp}$ denote $V$ shrunk and enlarged by $\eta_0$. The limit law of Step~3 of the proof of Proposition~\ref{prop:upper}, the fact that $Y$ has no atoms, and $\eta_0\to0$ give the limit for $T<\tau_m\le2T$; summing over dyadic ranges gives (b). Finally, $\mathbb P(|Y|>\pi)>0$ because $A(x_0)>\pi$ and the law of $Y$ has support $[-A(x_0),A(x_0)]$, by the argument of Lemma~\ref{lem:arcsin}(b).
\end{proof}

\begin{remark}\label{rem:upper}
The proof of Theorem~\ref{thm:order} needs the interpolating phase $\Upsilon$ to extend analytically to a strip. This holds for $c(-s)^{\kappa}$, but it does not follow from the real-variable hypotheses (H1$'$) and (H2), which is why Proposition~\ref{prop:upper} is stated separately.
\end{remark}

\begin{conjecture}[The boundary]\label{conj:boundary}
For $f(s)=c(-s)^{\kappa_1}$, Gram's law holds for all large $m$, whatever the value of $c\neq0$.
\end{conjecture}

\emph{Heuristic.} When $\kappa=\kappa_1$, a failure at height about $t$ requires $\arg\zeta(x_c+it)$ to lie within $O(\log\log t/\log t)$ of its supremum $A(x_c)$, where $x_c=\sigma_1+O(1/\log t)$. Each term of $A$ is attained only when $t\log p$ is close to a prescribed phase. The tail of $A$ beyond $y$ is $\asymp y^{1-\sigma_1}/\log y$. So this requires simultaneous approximation of $t\log p$ for every prime $p$ up to $(\log t)^{1/(\sigma_1-1)+o(1)}$, that is, up to $(\log t)^{29.5+o(1)}$. The probability of that at a given height is $\exp(-(\log t)^{29.5+o(1)})$. Summed over the $\asymp T\log T$ Gram intervals, the expected number of failures converges, so by Borel--Cantelli only finitely many are expected.

\section{Targets of order at most \texorpdfstring{$t^{1/2}$}{t\^{}(1/2)}: the regime \texorpdfstring{$\tfrac12<x_c\le1$}{1/2<xc<=1}}\label{sec:smallkappa}

Here $f(s)=c(-s)^{\kappa}$ with $c\ne0$ and $0<\kappa\le\tfrac12$, and we use the frame $(-1,2,\varphi_0)$ of Section~\ref{sec:counter}, case (ii), with $\varphi_0=\arg c-\kappa\pi/2$. The balance abscissa is now $x_c\approx\tfrac12+\kappa\in(\tfrac12,1]$, where the Euler product no longer converges absolutely. Both the single-crossing hypothesis of Proposition~\ref{prop:ident} and the control of $\arg\zeta$ near $x_c$ need new input. The difficulty is that $\zeta$ may have zeros between the critical line and $x_c$; Selberg's zero-density theorem shows that there are too few of them to matter at almost all separators.

\begin{theorem}[Small exponents]\label{thm:smallkappa}
Let $f(s)=c(-s)^{\kappa}$ with $c\neq0$ and $0<\kappa\le\tfrac12$. Put $x_0=\tfrac12+\kappa$ and $Y=Y_{x_0}$ \textup{(}Section~\textup{\ref{sec:main})}.
\begin{enumerate}
\item[(a)] For every $j\in\mathbb Z$,
\[
\lim_{T\to\infty}\frac{\#\{m:\ T<\tau_m\le2T,\ k_m=j\}}{\#\{m:\ T<\tau_m\le2T\}}=\mathbb P\bigl((2j-1)\pi<-Y<(2j+1)\pi\bigr).
\]
\item[(b)] $\#\{m:\ \tau_m\le T,\ n_m\ne1\}=o(T\log T)$.
\item[(c)] The law of $Y$ has support $\mathbb R$. Hence $\sup_m|k_m|=\infty$, and Gram's law fails for infinitely many $m$.
\end{enumerate}
\end{theorem}

Compared with Theorem~\ref{thm:order}, the counter now has unbounded support, but we only know that the failure count is $o(T\log T)$ and tends to infinity. Together with Theorems~\ref{thm:gram} and~\ref{thm:sharp}, for $f=c(-s)^{\kappa}$ Gram's law holds in every large interval when $\kappa>\kappa_1$, and fails infinitely often for every $0<\kappa<\kappa_1$.

\subsection{Tools}\label{sec:rh}

In this subsection and in Section~\ref{sec:gramlog}, $f$ is either $c(-s)^{\kappa}$ with $c\ne0$ and $0<\kappa\le\tfrac12$, or a target of Section~\ref{sec:slow} that satisfies \textup{(S3)}, with the frame (ii) or (iv) of Section~\ref{sec:counter}. Put $L_\tau=\log(\tau/2\pi)$. Constants $C$ may depend on $f$.

Recall that $S(t;x)=\frac1\pi\arg\zeta(x+it)$, defined by continuous variation from $2+it$. If $\zeta$ has no zeros in $[x_1,\infty)\times J$, for an interval $J$, then there $\pi S(t;x)=\Im\log\zeta(x+it)$ for the branch of $\log\zeta$ that is analytic there and agrees with $\sum_n\Lambda(n)n^{-s}/\log n$ for $x\ge2$. In particular $\pi\partial_tS=\Re\,\zeta'/\zeta$ and $\pi\partial_xS=\Im\,\zeta'/\zeta$ there.

\begin{lemma}\label{lem:hadamard}
For $w=x+iv$ with $\tfrac12<x\le3$ and $|v|\ge2$,
\begin{equation}\label{eq:hadU}
\Re\frac{\zeta'}{\zeta}(w)=\sum_\rho\frac{x-\beta}{|w-\rho|^2}-\frac12\log\frac{|v|}{2\pi}+O\Bigl(\frac1{|v|}\Bigr),\qquad\Bigl|\Bigl(\frac{\zeta'}{\zeta}\Bigr)'(w)\Bigr|\le\sum_\rho\frac1{|w-\rho|^2}+O\Bigl(\frac1{|v|}\Bigr).
\end{equation}
If RH holds, there is an absolute constant $C$ such that, for the same $w$,
\[
\Re\,\frac{\zeta'}{\zeta}(w)\ge-\frac12\log\frac{|v|}{2\pi}-\frac{C}{|v|}.
\]
\end{lemma}
\begin{proof}
By the Hadamard product \cite[\S2.12]{Titchmarsh},
\[
\frac{\zeta'}{\zeta}(w)=B+\sum_\rho\Bigl(\frac1{w-\rho}+\frac1\rho\Bigr)-\frac1{w-1}-\frac12\frac{\Gamma'}{\Gamma}\Bigl(\frac w2+1\Bigr)+\frac12\log\pi,
\]
with $\Re B=-\sum_\rho\Re(1/\rho)$. By Stirling's formula, $\Re\frac{\Gamma'}{\Gamma}(\frac w2+1)=\log\frac{|v|}2+O(1/|v|)$ and $(\frac{\Gamma'}{\Gamma})'(\frac w2+1)=O(1/|v|)$. Hence
\[
\Re\frac{\zeta'}{\zeta}(w)=\sum_\rho\Re\frac1{w-\rho}-\frac12\log\frac{|v|}{2\pi}+O\Bigl(\frac1{|v|}\Bigr),\qquad \Bigl(\frac{\zeta'}{\zeta}\Bigr)'(w)=-\sum_\rho\frac1{(w-\rho)^2}+O\Bigl(\frac1{|v|}\Bigr).
\]
Since $\Re\frac1{w-\rho}=\frac{x-\beta}{|w-\rho|^2}$, this gives \eqref{eq:hadU}. On RH, $\Re\frac1{w-\rho}=\frac{x-1/2}{|w-\rho|^2}>0$. This gives the last claim.
\end{proof}

\begin{lemma}[Identification]\label{lem:ident}
Let $f$ and its frame be as above. There are $m_1=m_1(f)$ and $C=C(f)$ such that the following holds for every $m\ge m_1$. Write $\tau=\tau_m$ and $R(\sigma)=|\zeta(\sigma+i\tau)/f(\sigma+i\tau)|$, and suppose that
\begin{equation}\label{eq:G1}
|\zeta(\sigma+i\tau)|\le\tfrac13|f(\sigma+i\tau)|\qquad(\tfrac12\le\sigma\le2).
\end{equation}
\begin{enumerate}
\item[(a)] For $\tfrac12<x\le2$,
\[
\mathcal R(x):=\log R(1-x)=(x-\tfrac12)L_\tau+\log|\zeta(x+i\tau)|-\log|f(1-x+i\tau)|+O(1/\tau),
\]
and, where $\zeta(x+i\tau)\ne0$, $\mathcal R'(x)=L_\tau+\Re\,\frac{\zeta'}{\zeta}(x+i\tau)+O(1/\tau)$.
\item[(b)] Suppose further that there is $\eta\in[0,\tfrac14]$ such that $\mathcal R<0$ on $(\tfrac12,\tfrac12+\eta]$, and that $\zeta(x+i\tau)\ne0$ and $\Re\,\frac{\zeta'}{\zeta}(x+i\tau)\ge-\tfrac23L_\tau$ for $\tfrac12+\eta<x\le2$. Then $\mathcal R'\ge\tfrac14L_\tau$ on $(\tfrac12+\eta,2]$, and $\mathcal R$ has a unique zero $x_c$ in $(\tfrac12,2]$; $\mathcal R<0$ on $(\tfrac12,x_c)$ and $\mathcal R>0$ on $(x_c,2]$. Moreover $k_m$ is the integer nearest to $-\tfrac12S(\tau;x_c)+e_m$, where $|e_m|\le\varepsilon_f(\tau)$, with $\varepsilon_f(\tau)\to0$ as $\tau\to\infty$ and $\varepsilon_f(\tau)\le C/\log\tau$ when $f=c(-s)^{\kappa}$ or $f=\ell$. In particular $|k_m+\tfrac12S(\tau;x_c)|\le\tfrac12+\varepsilon_f(\tau)$.
\item[(c)] On RH, the hypothesis of (b) holds with $\eta=0$.
\end{enumerate}
\end{lemma}
\begin{proof}
(a) Let $s=\sigma+i\tau$ with $-1\le\sigma<\tfrac12$, and write $\zeta(s)=\chi(s)\zeta(1-s)$. By Lemma~\ref{lem:chi} and its proof, $\log|\chi(s)|=(\tfrac12-\sigma)L_\tau+O(1/\tau)$ and $\partial_\sigma\log|\chi(s)|=\Re\,\chi'/\chi(s)=-L_\tau+O(1/\tau)$. Also $\partial_\sigma\log|f(s)|=\Re f'/f(s)=O(1/\tau)$, by \textup{(S2)} for the targets of Section~\ref{sec:slow}. Since $|\zeta(1-s)|=|\zeta(1-\sigma+i\tau)|$, putting $x=1-\sigma$ gives both formulas in (a).

(b) For $\tfrac12+\eta<x\le2$, (a) gives $\mathcal R'(x)\ge\tfrac13L_\tau-C/\tau\ge\tfrac14L_\tau$. By \eqref{eq:G1}, $\mathcal R(x)\to\log R(\tfrac12)\le\log\tfrac13$ as $x\to\tfrac12^+$. At $x=2$ we have $\mathcal R(2)=\tfrac32L_\tau+O(1)-\log|f(-1+i\tau)|\to\infty$, since $\log|f(-1+i\tau)|\le\tfrac12\log\tau+O(1)$ (by \textup{(S1)} for the targets of Section~\ref{sec:slow}). If $\eta=0$, then $\mathcal R$ is increasing on $(\tfrac12,2]$; if $\eta>0$, then $\mathcal R<0$ on $(\tfrac12,\tfrac12+\eta]$ and $\mathcal R$ is increasing on $[\tfrac12+\eta,2]$. Either way the claims about $x_c$ follow. Put $\sigma_c=1-x_c\in(-1,\tfrac12)$. Then $R>1$ on $[-1,\sigma_c)$ and, with \eqref{eq:G1}, $R<1$ on $(\sigma_c,2]$. This is the hypothesis of Proposition~\ref{prop:ident}. Put $\varepsilon_0(\tau)=\sup_{-1\le\sigma\le2}|\arg f(\sigma+i\tau)-\varphi_0|$. This tends to $0$ by \textup{(S3)}, and it is $O(1/\log\tau)$ for $f=c(-s)^{\kappa}$ and $f=\ell$. So the last claims follow from Proposition~\ref{prop:ident}, with $\varepsilon_f(\tau)=\varepsilon_0(\tau)/2\pi+C/\tau$.

(c) On RH, $\zeta$ has no zeros with $x>\tfrac12$, and Lemma~\ref{lem:hadamard} gives $\Re\,\frac{\zeta'}{\zeta}(x+i\tau)\ge-\tfrac12L_\tau-C/\tau\ge-\tfrac23L_\tau$.
\end{proof}

On RH, parts (b) and (c) show that $R$ is strictly decreasing on $[-1,\tfrac12)$ whenever \eqref{eq:G1} holds. This is a quantitative form, for the ratio $|\zeta/f|$, of the horizontal monotonicity of $|\zeta|$: Matiyasevich, Saidak and Zvengrowski \cite{MSZ} proved that $|\zeta(\sigma+it)|$ is strictly decreasing in $\sigma$ for $\sigma\le0$ and $|t|\ge8$, and that extending this to $\sigma<\tfrac12$ is equivalent to RH. Without RH, we verify the hypothesis of (b) at almost all separators, using Selberg's zero-density theorem; see the proofs of Theorems~\ref{thm:smallkappa} and~\ref{thm:gramfail}.

\subsection{Proof of Theorem~\ref{thm:smallkappa}}

\begin{proof}[Proof of Theorem~\ref{thm:smallkappa}]
Put $c^*=\log|c|+\kappa\log2\pi$ and $L=\log T$, with $T$ large. Fix $B\ge1$, put $A=|c^*|+B+3$, and fix $\epsilon\in(0,1/(100A)]$. Let $h=4\pi/L$ and $I=[x_0-2A/L,\,x_0+2A/L]$. For real $t$ put $Q(t)=I\times[t-h,t+h]$ and $F(t)=\max_{Q(t)}|\zeta'/\zeta|$. Consecutive separators near height $T$ are at least $2\pi/(L+1)$ and at most $g:=2\pi/(L-2)$ apart; we call the set of $t$ whose nearest separator is $\tau_m$ the \emph{cell} of $\tau_m$. Note that $\tfrac32g\le h$. Let $\mathcal Z$ be the set of zeros $\rho=\beta+i\gamma$ of $\zeta$ with $\beta>\tfrac12+\kappa/4$.

\emph{Step 1: limit laws.} By the Bohr--Jessen theorem \cite{BohrJessen} (see also \cite{Laurincikas}), there is a probability measure $\mu$ on $\mathbb C$, namely the law of $-\sum_p\log(1-p^{-x_0}e^{i\vartheta_p})$, with
\[
\frac1T\operatorname{meas}\{t\in[T,2T]:\ \log\zeta(x_0+it)\in R\}\to\mu(R)
\]
for every rectangle $R$ with sides parallel to the axes. Here $\log\zeta(x_0+it)=\log|\zeta(x_0+it)|+i\pi S(t;x_0)$. (The proof in \cite{BohrJessen}, stated for $t\in[-T,T]$, applies equally to $t\in[T,2T]$.) Since $\tfrac12<x_0\le1$, the support of $\mu$ is $\mathbb C$ \cite{BohrJessen,Laurincikas}. The imaginary part of a $\mu$-distributed variable is $Y$, so the law of $Y$ has support $\mathbb R$. Since $\mu$ is tight, the same limit relation holds for every set $\{z:\ \Im z\in U\}$ with $U$ a finite union of intervals. Finally, $Y$ has no atoms: $Y=Y_2+Y'$ with $Y_2=-\arg(1-2^{-x_0}e^{i\vartheta_2})$ independent of $Y'$, and $Y_2$ has no atoms.

By Guo's theorem \cite{Guo} (see also \cite{LesterLD}), which is unconditional and applies at the fixed abscissa $x_0$, $\zeta'/\zeta(x_0+it)$ has a limiting distribution: there is a probability measure $\nu$ on $\mathbb C$ with $\frac1T\int_0^TQ\bigl(\frac{\zeta'}{\zeta}(x_0+it)\bigr)\,dt\to\int Q\,d\nu$ for every smooth, compactly supported $Q$. Taking $0\le Q\le1$ with $Q=1$ on $|z|\le B/2$ and $Q=0$ for $|z|\ge B$, and applying this with $T$ and $2T$, we find that the set of $t\in[T,2T]$ with $|\zeta'/\zeta(x_0+it)|>B$ has measure at most $(2\nu(\{|z|>B/2\})+o(1))T$.

\emph{Step 2: the exceptional set.} Let $E_T$ be the set of $t\in[T,2T]$ such that
\begin{enumerate}
\item[(a)] $\log\zeta(x_0+it)\in[-B,B]^2$;
\item[(b)] $\max\{|\zeta(\sigma+iv)|:\ \tfrac12\le\sigma\le2,\ |v-t|\le h\}\le T^{\kappa/2}$;
\item[(c)] $|\zeta'/\zeta(x_0+it)|\le B$;
\item[(d)] no $\rho\in\mathcal Z$ has $|\gamma-t|\le1$;
\item[(e)] $\widetilde N(t):=4\sum_{\rho\in\mathcal Z,\,|\gamma-t|>1}(\gamma-t)^{-2}\le1$.
\end{enumerate}
We claim that $\operatorname{meas}([T,2T]\setminus E_T)\le(\beta(B)+o(1))T$, where $\beta(B)=1-\mu([-B,B]^2)+2\nu(\{|z|>B/2\})\to0$ as $B\to\infty$. For (a) and (c) this is Step~1. For (b), note that (b) holds whenever $m_T(t)\le\tfrac12\kappa\log T$, so by Lemma~\ref{lem:max} the exceptional set has measure $O(T\log\log T/(\kappa\log T)^2)$. For (d) and (e), let $\mathcal Z'$ be the set of $\rho\in\mathcal Z$ with $0<\gamma\le3T$. Selberg's density theorem \cite[Thm.~1]{Selberg46} gives $\#\mathcal Z'\ll T^{1-\kappa/16}\log T$. The set where (d) fails has measure at most $2\,\#\mathcal Z'$. The zeros of $\mathcal Z$ outside $\mathcal Z'$ contribute $O(T^{-1}\log T)$ to $\widetilde N$, and $\int_T^{2T}\sum_{\rho\in\mathcal Z',|\gamma-t|>1}(\gamma-t)^{-2}\,dt\le2\,\#\mathcal Z'$, so (e) fails on a set of measure $O(T^{1-\kappa/16}\log T)$.

\emph{Step 3: a derivative bound.} Let $t\in E_T$. By (d), $\zeta$ has no zeros $x+iv$ with $x>\tfrac12+\kappa/4$ and $|v-t|\le1$. Let $z=x+iv$ with $x>\tfrac12+\kappa/4$ and $|v-t|\le h$. Every zero with $\beta>x$ lies in $\mathcal Z$ and, by (d), has $|\gamma-v|\ge\tfrac12|\gamma-t|$. So, by (e) and \eqref{eq:hadU},
\begin{equation}\label{eq:negpart2}
\sum_{\beta>x}\frac{\beta-x}{|z-\rho|^2}\le\widetilde N(t)\le1,\qquad \Re\frac{\zeta'}{\zeta}(z)\ge-\frac12\log\frac{|v|}{2\pi}-2 .
\end{equation}
We claim that $F(t)\le K:=2B+(8A+16\pi)/\kappa+2$ for large $T$. The box $Q(t)$ contains no zeros, so $F(t)<\infty$. Let $z\in Q(t)$; then $x\ge\tfrac12+\tfrac34\kappa$ for large $T$. In the second sum in \eqref{eq:hadU}, the zeros with $\beta>x-\kappa/2$ lie in $\mathcal Z$ and contribute at most $\widetilde N(t)\le1$. The zeros with $\beta\le x-\kappa/2$ satisfy $|z-\rho|^{-2}\le(2/\kappa)(x-\beta)|z-\rho|^{-2}$, and by the first identity in \eqref{eq:hadU} and \eqref{eq:negpart2} their contribution is at most
\[
\frac2\kappa\sum_{\beta<x}\frac{x-\beta}{|z-\rho|^2}\le\frac2\kappa\Bigl(F(t)+\frac12\log\frac{v}{2\pi}+1+O\Bigl(\frac1T\Bigr)\Bigr)\le\frac2\kappa\bigl(F(t)+L\bigr).
\]
So $|(\zeta'/\zeta)'(z)|\le(2/\kappa)(F(t)+L)+2$. The segment from $x_0+it$ to $z$ lies in $Q(t)$ and has length at most $(2A+4\pi)/L$, so by (c),
\[
F(t)\le B+\frac{2A+4\pi}{L}\Bigl(\frac2\kappa\bigl(F(t)+L\bigr)+2\Bigr)\le B+\frac{4A+8\pi}{\kappa}+1+\frac{4A+8\pi}{\kappa L}F(t),
\]
which gives the claim. In particular $F(t)\le\epsilon L$ once $T\ge T_1(\epsilon,A,B,\kappa)$. Since $\zeta$ has no zeros in $[x_0,\infty)\times[t-h,t+h]$, it follows that for $|t'-t|\le h$ we have $|\log\zeta(x_0+it')-\log\zeta(x_0+it)|\le hF(t)\le4\pi\epsilon<0.2$, and in particular $|S(t';x_0)-S(t;x_0)|\le4\epsilon$.

\emph{Step 4: good separators.} Call $\tau_m$ \emph{good} if its cell meets $E_T$. Let $t\in E_T$, and let $\tau=\tau_m$ be the separator nearest to $t$, so that $|\tau-t|\le h$.
\begin{enumerate}
\item[(i)] For $\tfrac12\le\sigma\le2$, $|\zeta(\sigma+i\tau)|\le T^{\kappa/2}\le\tfrac13|c|(T-1)^\kappa\le\tfrac13|f(\sigma+i\tau)|$. So \eqref{eq:G1} holds.
\item[(ii)] By Step~3, $|\log|\zeta(x_0+i\tau)||\le B+0.2$.
\item[(iii)] We verify the hypothesis of Lemma~\ref{lem:ident}(b) with $\eta=\kappa/4$. Since $\log|f(1-x+i\tau)|=\log|c|+\kappa\log\tau+O(\tau^{-2})$, Lemma~\ref{lem:ident}(a) and condition (b) of Step~2 give, for $\tfrac12<x\le\tfrac12+\kappa/4$,
\[
\mathcal R(x)\le\tfrac14\kappa L_\tau+\tfrac12\kappa L-\kappa\log\tau-\log|c|+o(1)\le-\tfrac14\kappa L+O(1)<0 .
\]
For $\tfrac12+\kappa/4<x\le2$, $\zeta(x+i\tau)\ne0$ and, by \eqref{eq:negpart2}, $\Re\,\zeta'/\zeta(x+i\tau)\ge-\tfrac12L_\tau-2\ge-\tfrac23L_\tau$. So Lemma~\ref{lem:ident}(b) applies: $\mathcal R$ is increasing on $[\tfrac12+\kappa/4,2]$, with a unique zero $x_c$, and $k_m$ is the integer nearest to $-\tfrac12S(\tau;x_c)+e_m$ with $|e_m|\le C/L$.
\item[(iv)] Since $(x-\tfrac12)L_\tau-\kappa\log\tau=(x-x_0)L_\tau-\kappa\log2\pi$, Lemma~\ref{lem:ident}(a) gives
\[
\mathcal R(x)=(x-x_0)L_\tau-c^*+\log|\zeta(x+i\tau)|+O(1/\tau).
\]
For $x\in I$ we have $|\log|\zeta(x+i\tau)|-\log|\zeta(x_0+i\tau)||\le(2A/L)F(t)\le2A\epsilon\le0.02$. Since $L_\tau\ge L-2$, at $x=x_0\pm2A/L$ we get $\pm\mathcal R(x)\ge2A(1-2/L)-|c^*|-B-0.22-o(1)>0$. As $I\subset[\tfrac12+\kappa/4,2]$, it follows that $x_c\in I$.
\item[(v)] $|S(\tau;x_c)-S(\tau;x_0)|\le(2A/\pi L)F(t)\le A\epsilon$.
\item[(vi)] By (iii) and (v), $k_m$ is the integer nearest to $-\tfrac12S(\tau;x_0)+e_m'$, with $|e_m'|\le\tfrac12A\epsilon+C/L\le A\epsilon$.
\end{enumerate}
Each cell has length at least $2\pi/(L+1)$, and the cell of a separator that is not good lies in $[T-1,2T+1]\setminus E_T$. Hence the number of $m$ with $T<\tau_m\le2T$ for which $\tau_m$ is not good is at most $(\beta(B)+o(1))TL/2\pi+O(1)$.

\emph{Step 5: proof of (a).} Let $j\in\mathbb Z$ and $\eta_0\ge10\pi A\epsilon$. If $t\in E_T$ and $-\pi S(t;x_0)\in[(2j-1)\pi+\eta_0,(2j+1)\pi-\eta_0]$, then the separator $\tau_m$ nearest to $t$ satisfies $-\pi S(\tau_m;x_0)\in[(2j-1)\pi+\eta_0-4\pi\epsilon,(2j+1)\pi-\eta_0+4\pi\epsilon]$ by Step~3, and hence $k_m=j$ by Step~4(vi). Each separator is nearest to a set of $t$ of measure at most $g$. So, by Steps~1 and~2,
\[
\liminf_{T\to\infty}\frac{\#\{m:\ T<\tau_m\le2T,\ k_m=j\}}{TL/2\pi}\ge\mathbb P\bigl(-Y\in[(2j-1)\pi+\eta_0,(2j+1)\pi-\eta_0]\bigr)-\beta(B).
\]
Letting $\epsilon\to0$, then $\eta_0\to0$, then $B\to\infty$, the right-hand side tends to $p_j:=\mathbb P((2j-1)\pi<-Y<(2j+1)\pi)$. Since $Y$ has no atoms, $\sum_jp_j=1$. As $\#\{m:T<\tau_m\le2T\}=(1+o(1))TL/2\pi$, these lower bounds for all $j$ force the limits in (a).

\emph{Step 6: proof of (b).} Let $T<\tau_m\le2T$ with $n_m\ne1$, so that $k_{m-1}\ne k_m$. At most $2(\beta(B)+o(1))TL/2\pi+O(1)$ such $m$ have $\tau_{m-1}$ or $\tau_m$ not good. Suppose both are good, and let $t\in E_T$ lie in the cell of $\tau_m$. Then $[\tau_{m-1},\tau_m]\subset[t-\tfrac32g,t+\tfrac12g]\subset[t-h,t+h]$, so by Step~3, $|S(t';x_0)-S(t;x_0)|\le4\epsilon$ for all $t'\in[\tau_{m-1},\tau_m]$. By Step~4(vi), a point of $\tfrac12+\mathbb Z$ lies between $-\tfrac12S(\tau_{m-1};x_0)+e_{m-1}'$ and $-\tfrac12S(\tau_m;x_0)+e_m'$. It follows that $\operatorname{dist}(S(t';x_0),2\mathbb Z+1)\le12A\epsilon$ for all $t'\in[\tau_{m-1},\tau_m]$. These intervals are disjoint and have length at least $2\pi/(L+1)$, so by Step~1
\[
\limsup_{T\to\infty}\frac{\#\{m:\ T<\tau_m\le2T,\ n_m\ne1\}}{TL/2\pi}\le2\beta(B)+\mathbb P\bigl(\operatorname{dist}(Y/\pi,2\mathbb Z+1)\le12A\epsilon\bigr).
\]
Let $\epsilon\to0$ (using again that $Y$ has no atoms), then $B\to\infty$, and sum over dyadic ranges.

\emph{Step 7: proof of (c).} Since the law of $Y$ has support $\mathbb R$, (a) shows that for every $J$ a positive proportion of separators have $|k_m|>J$. If $n_m=1$ for all $m\ge m_2$, then $k_m=k_{m_2}$ for all $m\ge m_2$, which is a contradiction.
\end{proof}
\section{Slowly growing targets: counting and horizontal distribution}\label{sec:slow}

\subsection*{Hypotheses}
Throughout this section, $f$ is analytic on $\{t>t_0\}$ with $|f(s)|\le C_1|s|^{K}$ there. Moreover, as $t\to\infty$,
\begin{enumerate}
\item[(S1)] $|f(\sigma+it)|\to\infty$ uniformly for $\sigma\ge-1$, and $\log|f(\sigma+it)|=o(\log t)$ uniformly for $\sigma$ in compact sets;
\item[(S2)] $f'(\sigma+it)\ll|f(\sigma+it)|/t$ uniformly for $\sigma$ in compact sets.
\end{enumerate}
Put $L_f(t)=\log|f(it)|$. By (S2), $f(\sigma+it)=f(it)\bigl(1+O(1/t)\bigr)$ for bounded $\sigma$, and $L_f(2t)-L_f(t)=O(1)$. Examples include
\[
\ell(s)=\log(-s),\qquad c\,\ell(s)^{A}\ (A>0),\qquad \log\ell(s),\qquad \exp\bigl(\sqrt{\ell(s)}\bigr),
\]
with $L_f(t)$ asymptotic to $\log\log t$, $A\log\log t$, $\log\log\log t$ and $\sqrt{\log t}$ respectively. For the results on Gram's law in Section~\ref{sec:gramlog} we also assume
\begin{enumerate}
\item[(S3)] $\arg f(\sigma+it)\to\varphi_0$ as $t\to\infty$, uniformly for $-1\le\sigma\le2$, for a continuous branch of $\arg f$.
\end{enumerate}
This holds for all four examples, with $\varphi_0=\arg c$ for $c\,\ell(s)^A$ and $\varphi_0=0$ for the others.

Write $H_f=\zeta-f$, and let $\rho=\beta+i\gamma$ denote its zeros with $\gamma>\max(t_0,0)$ and $-1<\beta<2$; sums over $\rho$ run over these zeros, and $N_f(T)$ is their number with $\gamma\le T$. By Lemma~\ref{lem:strip} below, for large $\gamma$ every zero of $H_f$ lies in this strip. For the principal example we write $H=H_\ell$ and $N_H=N_\ell$. By \eqref{eq:modulus},
\begin{equation}\label{eq:betaH}
\beta=\tfrac12-\frac{\log|f(\rho)|-\log|\zeta(1-\rho)|}{\log(\gamma/2\pi)}+O(\gamma^{-1}).
\end{equation}

\begin{lemma}\label{lem:strip}
For every $\varepsilon>0$, all zeros of $H_f$ with $\gamma\ge t_0(f,\varepsilon)$ satisfy $-\varepsilon<\beta<1+\varepsilon$. If $K<\tfrac32$, there are no zeros with $\beta\le-1$ and $\gamma\ge t_0$; this holds in particular for $f=\ell$.
\end{lemma}
\begin{proof}
For $\sigma\ge1+\varepsilon$: $|\zeta|\le\zeta(1+\varepsilon)<|f|$ for large $t$, by (S1). Fix $\sigma_L<\min(-1,\tfrac12-K)$. For $\sigma_L\le\sigma\le-\varepsilon$: $|\zeta(s)|\gg_\varepsilon t^{1/2+\varepsilon}$ by Lemma~\ref{lem:chi}, while $|f|=t^{o(1)}$ by (S1). For $\sigma<\sigma_L$: use Lemma~\ref{lem:left}. If $K<\tfrac32$ we may take $\sigma_L=-1$.
\end{proof}

\begin{lemma}[Horizontal segments]\label{lem:horiz}
Fix $\sigma_0<0<2=c$. For $T$ large and not the ordinate of a zero of $H_f$, define $\arg H_f(\sigma+iT)$ by continuous variation from $c+iT$. Then $\arg H_f(\sigma+iT)=O(\log T)$ uniformly for $\sigma_0\le\sigma\le c$.
\end{lemma}
\begin{proof}
This is Backlund's argument (cf. \cite[\S9.4]{Titchmarsh}). Let $\omega=\arg H_f(c+iT)$. Let $k$ be the number of zeros of $\Re\bigl(e^{-i\omega}H_f(\sigma+iT)\bigr)$ on $[\sigma_0,c]$. Then $|\arg H_f(\sigma+iT)-\omega|\le\pi(k+1)$. Apply Jensen's formula to
\[
g(z)=\tfrac12\bigl(e^{-i\omega}H_f(z+iT)+\overline{e^{-i\omega}H_f(\bar z+iT)}\bigr)
\]
in a disc centred at $c$ of radius $c-\sigma_0+1$. There $|g|\ll T^{A}$, while $g(c)=|H_f(c+iT)|\ge|f(c+iT)|-\zeta(2)\gg1$. Hence $k\ll\log T$. Finally, $\omega=\arg(1-f(c+iT))+O(1)$ by Lemma~\ref{lem:right}, and $\arg(1-f(c+it))$ varies by $\ll\int^{T}|f'/f|\,dt\ll\log T$, by (S2).
\end{proof}

\begin{lemma}[Right edge]\label{lem:right}
Let $c=2$ and $u(t)=(\zeta(c+it)-1)/(1-f(c+it))$. Then $|u(t)|<1$ for large $t$, and
\[
\int_{T_0}^{T}\log|H_f(c+it)|\,dt=\int_{T_0}^{T}\log|1-f(it)|\,dt+O(\log T).
\]
For $f=\ell$ the error term is $C+O(1/\log T)$ for a constant $C$.
\end{lemma}
\begin{proof}
We have $|\zeta(2+it)-1|\le\zeta(2)-1<0.65$, and $|1-f(c+it)|\to\infty$. Hence $|u|<1$, and
\[
\log|H_f|=\log|1-f|+\Re\sum_{k\ge1}\frac{(-1)^{k+1}}{k}u^k.
\]
Now $(\zeta(c+it)-1)^k=\sum_{n\ge2^k}a_k(n)n^{-c-it}$, with $a_k(n)\ge0$ and $\sum_na_k(n)n^{-c}=(\zeta(c)-1)^k$. Put $g_k(t)=(1-f(c+it))^{-k}$ and $M=\min_{[T_0,T]}|1-f(c+it)|$. By (S2), $|g_k'|\ll k\,t^{-1}M^{-k}$, so integrating by parts gives $\int_{T_0}^{T}n^{-it}g_k(t)\,dt\ll kM^{-k}\log T/\log n$. Summing over $n$ and $k$ gives $O(\log T)$, because $(\zeta(c)-1)/M<1$. For $f=\ell$ we have $|g_k'|\ll k\,t^{-1}(\log t)^{-k-1}$ and $|g_k(T)|\ll(\log T)^{-k}$, so this part is a constant plus $O(1/\log T)$. Finally, $f(c+it)=f(it)(1+O(1/t))$ gives $\log|1-f(c+it)|=\log|1-f(it)|+O(1/t)$. For $f=\ell$ the difference is $\Re\bigl(2i/(t(1-\log(-it)))\bigr)+O(t^{-2})=O(1/(t\log^2t))$, which is integrable with tail $O(1/\log T)$.
\end{proof}

\begin{theorem}[Counting]\label{thm:countH}
$N_f(T)=\dfrac{\theta(T)}{\pi}+O(\log T)=\dfrac{T}{2\pi}\log\dfrac{T}{2\pi e}+O(\log T)$.
\end{theorem}
\begin{proof}
Apply the argument principle on $[-1,2]\times[T_0,T]$. By Lemma~\ref{lem:strip} no zeros lie on or near the vertical sides for large $t$.
\begin{itemize}
\item On $\sigma=-1$: $H_f=\zeta(1-f/\zeta)$ with $|f/\zeta|\ll t^{-3/2+o(1)}$, and $\arg\zeta=\arg\chi+\arg\zeta(2-it)$ with $|\arg\zeta(2-it)|\le A(2)$. Going downward, Lemma~\ref{lem:chi} gives the change $2\theta(T)-2\theta(T_0)+O(1)$.
\item On $\sigma=2$: $H_f=(1-f)(1+u)$ with $|u|<1$, so the change is $O(\log T)$ by (S2).
\item On the horizontal sides: $O(\log T)$ (Lemma~\ref{lem:horiz}).
\end{itemize}
\end{proof}

\begin{theorem}[Leftward drift]\label{thm:drift}
\[
\sum_{0<\gamma\le T}\bigl(\beta-\tfrac12\bigr)=-\frac1{2\pi}\int_{T_0}^{T}\log\bigl|1-f(it)\bigr|\,dt+O(\log T)=-\frac{1+o(1)}{2\pi}\,T\,L_f(T).
\]
In particular, the mean of $\beta-\tfrac12$ over $T<\gamma\le2T$ is $-(1+o(1))\,L_f(T)/\log T$. For $f=\ell$,
\[
\sum_{0<\gamma\le T}\bigl(\beta-\tfrac12\bigr)=-\frac1{2\pi}\int_{2}^{T}\log\bigl|1-\log(-it)\bigr|\,dt+O(\log T)=-\frac{T}{2\pi}\log\log T+O\Bigl(\frac{T}{\log T}\Bigr).
\]
\end{theorem}
\begin{proof}
Apply Littlewood's lemma \cite[\S9.9]{Titchmarsh} on $[\sigma_0,c]\times[T_0,T]$ with $\sigma_0=-1$ and $c=2$:
\begin{multline*}
2\pi\!\!\sum_{T_0<\gamma\le T}\!\!(\beta-\sigma_0)=\int_{T_0}^{T}\!\log|H_f(\sigma_0+it)|\,dt-\int_{T_0}^{T}\!\log|H_f(c+it)|\,dt\\
+\int_{\sigma_0}^{c}\!\bigl(\arg H_f(\sigma+iT)-\arg H_f(\sigma+iT_0)\bigr)\,d\sigma.
\end{multline*}
The last term is $O(\log T)$ by Lemma~\ref{lem:horiz}.

\emph{Left edge.} We have $\log|H_f(-1+it)|=\log|\chi(-1+it)|+\log|\zeta(2-it)|+O(t^{-3/2+o(1)})$. By Lemma~\ref{lem:chi}, $\log|\chi(-1+it)|=\tfrac32\log(t/2\pi)+O(1/t)$. Also
\[
\int_{T_0}^{T}\log|\zeta(2-it)|\,dt=\Re\sum_n\frac{\Lambda(n)}{n^2\log n}\int_{T_0}^{T}n^{it}\,dt=O(1).
\]
So this edge contributes $\tfrac32\int_{T_0}^{T}\log(t/2\pi)\,dt+O(\log T)$.

\emph{Right edge.} Use Lemma~\ref{lem:right}.

\emph{Conclusion.} Subtract $\tfrac32\cdot2\pi N_f(T)$ from both sides. By Theorem~\ref{thm:countH}, together with $2\theta(T)-2\theta(T_0)=\int_{T_0}^{T}\log(t/2\pi)\,dt+O(1)$, the $\log(t/2\pi)$ terms cancel up to $O(\log T)$. For the asymptotic form: by (S1), $\log|1-f(it)|=L_f(t)+O(1/|f(it)|)$. Since $L_f(2t)-L_f(t)=O(1)$ and $L_f\to\infty$, we have $\int_{T}^{2T}L_f=(1+o(1))\,T L_f(T)$. For $f=\ell$, use $\log|1-\log(-it)|=\log\log t-1/\log t+O(1/\log^2t)$.
\end{proof}

\begin{theorem}[Few solutions to the right]\label{thm:right}
For large $T$,
\[
\#\{T<\gamma\le2T:\ \beta>\tfrac12\}\ll\frac{T\log T\,\log\log T}{L_f(T)^2}+\log^2T .
\]
In particular, if $L_f(T)/\sqrt{\log\log T}\to\infty$, the proportion of the solutions with $T<\gamma\le2T$ that lie to the right of the critical line tends to $0$. For $f=\ell$ this proportion is $\ll1/\log\log T$.
\end{theorem}
\begin{proof}
Let $L=\log T$, $\lambda=\log\log T$, $w=L_f(T)/(2L)$ and $u=\tfrac12-w$. By (S1), $w\to0$. Each zero with $\beta>\tfrac12$ has $\beta-u>w$, so
\begin{equation}\label{eq:rightcount}
\#\{T<\gamma\le2T:\ \beta>\tfrac12\}\le w^{-1}\sum_{\substack{T<\gamma\le2T\\ \beta>u}}(\beta-u).
\end{equation}
We may assume that $T$ and $2T$ are not ordinates of zeros of $H_f$. Apply Littlewood's lemma on $[u,2]\times[T,2T]$. The horizontal sides contribute $O(\log T)$ by Lemma~\ref{lem:horiz}, and the right side is handled by Lemma~\ref{lem:right}. On the left side,
\[
\log|H_f(u+it)|\le\log|f(u+it)|+\log\Bigl(1+\frac{|\zeta(u+it)|}{|f(u+it)|}\Bigr).
\]
By (S2), $\log|f(u+it)|-\log|1-f(it)|=O(1/|f(it)|)+O(1/t)$, and $|f(it)|\asymp|f(iT)|$ for $T\le t\le2T$. So
\[
2\pi\sum_{\substack{T<\gamma\le2T\\ \beta>u}}(\beta-u)\le\int_T^{2T}\log\Bigl(1+\frac{|\zeta(u+it)|}{|f(u+it)|}\Bigr)dt+O\Bigl(\frac{T}{|f(iT)|}+\log T\Bigr).
\]
By the functional equation and Lemma~\ref{lem:chi}, $|\zeta(u+it)|=|\chi(u+it)||\zeta(1-u+it)|\le e^{wL+O(1/t)}|\zeta(\tfrac12+w+it)|$, and $wL=\tfrac12L_f(T)$. Also $\log|f(u+it)|\ge L_f(T)-C$. Put $a=\tfrac12L_f(T)$ and $X(t)=\log|\zeta(\tfrac12+w+it)|$. Then the integrand is at most $\log(1+Ce^{X-a})$. For $X\le a$ we bound this by $Ce^{X-a}$, and for $X>a\ge1$ by $(1+\log(1+C))X$.

By Lemma~\ref{lem:moment}, $\int_T^{2T}(X^+)^2\,dt\ll T\lambda$. So $\int_{X>a}X\,dt\le a^{-1}\int_T^{2T}(X^+)^2\,dt\ll T\lambda/a$, and, by Chebyshev's inequality, $\operatorname{meas}\{t\in[T,2T]:X(t)>j\}\ll T\lambda/j^2$ for $j\ge1$. Hence
\begin{align*}
\int_{X\le a}e^{X-a}\,dt&\le e^{1-a}T+e^{-a}\sum_{2\le j\le a+1}e^{j}\operatorname{meas}\{X>j-1\}\\
&\ll Te^{-a}+T\lambda e^{-a}\sum_{2\le j\le a+1}\frac{e^{j}}{(j-1)^2}\ll Te^{-a}+\frac{T\lambda}{a^2}.
\end{align*}
Altogether $\sum_{\beta>u}(\beta-u)\ll T\lambda/L_f(T)+\log T$, since $1/|f(iT)|=e^{-L_f(T)}\le\lambda/L_f(T)$ for large $T$. By \eqref{eq:rightcount} the number of zeros in question is $\ll(L/L_f(T))(T\lambda/L_f(T)+\log T)$, which gives the first claim. The proportion follows from Theorem~\ref{thm:countH}, and for $f=\ell$ we have $L_f(T)=\log\log T+O(1/\log^2T)$.
\end{proof}

\begin{remark}\label{rem:split}
For fixed $a\neq0$, Selberg's conjecture (see \cite{Lester}) puts about a quarter of the $a$-points to the right of $\tfrac12$, at distances of order $\sqrt{\log\log T}/\log T$. Theorem~\ref{thm:right} shows that for a growing target the proportion to the right tends to $0$ as soon as $L_f(T)/\sqrt{\log\log T}\to\infty$. The condition is natural. A solution with $\beta>\tfrac12$ is a point to the right of the critical line where $\log|\zeta|=\log|f|\approx L_f(T)$, and by Selberg's theorem $\log|\zeta|$ has typical size $\sqrt{\tfrac12\log\log T}$ there. Accordingly, the proof puts the line of Littlewood's lemma at $\tfrac12-L_f(T)/(2\log T)$, half way to the mean abscissa $\tfrac12-L_f(T)/\log T$ of the solutions (Theorem~\ref{thm:drift}). When $L_f(T)\ll\sqrt{\log\log T}$ we expect a positive proportion of the solutions to lie to the right of the line (Section~\ref{sec:open}).
\end{remark}

\section{Slowly growing targets: Gram's law}\label{sec:gramlog}

In this section $f$ is a target of Section~\ref{sec:slow} that satisfies \textup{(S3)}, and we use the frame $(-1,2,\varphi_0)$ of Section~\ref{sec:counter}, case (iv). So $\theta(\tau_m)=(m+\tfrac12)\pi-\tfrac12\varphi_0$, and $n_m$ counts the zeros of $H_f$ with $\tau_{m-1}<t\le\tau_m$. Put $S_f(T)=N_f(T)-(\theta(T)+\tfrac12\varphi_0)/\pi$. By Proposition~\ref{prop:counter} and Corollary~\ref{cor:block},
\begin{equation}\label{eq:SHk}
S_f(\tau_m)=C_f-\tfrac12-k_m ,
\end{equation}
and $S_f(T)=O(\log T)$ by Theorem~\ref{thm:countH}. For $f=\ell$ we have $\varphi_0=0$, and we write $H=\zeta-\ell$, $N_H=N_\ell$, $S_H=S_\ell$ and $C_H=C_\ell$, as in Theorem~D.

Since $N(T)=\theta(T)/\pi+1+S(T)$, \eqref{eq:SHk} can be rewritten as
\begin{equation}\label{eq:twoerr}
N_f(\tau_m)-N(\tau_m)=-S(\tau_m)-k_m+C_f',\qquad C_f'=C_f-\tfrac32+\frac{\varphi_0}{2\pi}.
\end{equation}
So at the separators the solutions follow $\theta(T)/\pi$ as far as the counter $k_m$ stays bounded, and they follow the zeros of $\zeta$ as far as $k_m$ cancels $S(\tau_m)$. They cannot follow both, whatever $f$ is: the two counting errors differ by $1+S(\tau_m)$ up to a constant, and $S$ is unbounded. The content of Theorem~\ref{thm:gramfail} is that, once $L_f(T)/\sqrt{\log\log T}\to\infty$, the solutions do not follow $\theta(T)/\pi$; for $f=\ell$, Corollary~\ref{cor:neither} shows that they do not follow the zeros of $\zeta$ either. The reason is visible in \eqref{eq:twoerr}: at almost all separators $k_m=-\tfrac12S(\tau_m;x_c)+O(1)$ (Step~3 of the proof of Theorem~\ref{thm:gramfail}), so that
\[
N_f(\tau_m)-N(\tau_m)=-S(\tau_m)+\tfrac12S(\tau_m;x_c)+O(1).
\]

\subsection{The central limit theorem}

\begin{theorem}[Distribution of the counting error]\label{thm:gramfail}
Suppose that $L_f(T)/\sqrt{\log\log T}\to\infty$. Let $\delta_T>0$ satisfy $\delta_T\log T=L_f(T)+O(1)$, and put
\[
\Psi_T=\frac12\sum_{p<T}p^{-1-2\delta_T}=\frac12\log\frac{\log T}{L_f(T)}+O(1),\qquad \sigma_T=\frac{\sqrt{\Psi_T}}{2\pi}.
\]
Then for all fixed $a<b$,
\[
\frac{\#\{m:\ T<\tau_m\le2T,\ S_f(\tau_m)\in[a\sigma_T,b\sigma_T]\}}{\#\{m:\ T<\tau_m\le2T\}}\longrightarrow\frac1{\sqrt{2\pi}}\int_a^be^{-x^2/2}\,dx .
\]
In particular $|S_f(\tau_m)|\to\infty$ in density, $S_f$ is unbounded, and Gram's law for $\zeta(s)=f(s)$ fails for infinitely many $m$. For $f=\ell$ one may take $\delta_T=\log\log T/\log T$, and then $\Psi_T=\frac12\sum_{p<T}p^{-1-2\log\log T/\log T}=\frac12\log\log T+O(\log\log\log T)$.
\end{theorem}

For $f=\ell$ the variance is a quarter of the variance $\log\log T/2\pi^2$ of $S(t)$ in Selberg's theorem, in line with $S_H\approx\tfrac12S(\cdot\,;x_c)$ (Remark~\ref{rem:SH}).

The proof identifies $k_m$ with $-\tfrac12S(\tau_m;x_c)$ as in Lemma~\ref{lem:ident}, where $x_c=\tfrac12+(1+o(1))\delta_T$. It then uses the Selberg--Tsang central limit theorem near the critical line \cite{Tsang} and Lester's limit theorem for $\zeta'/\zeta$ \cite{LesterLD}, both of which are unconditional, and Selberg's zero-density theorem in place of RH. Lemma~\ref{lem:max} keeps $|\zeta|$ below $|f|$ on the right half of almost every separator, and this is where the hypothesis $L_f(T)/\sqrt{\log\log T}\to\infty$ enters (Remark~\ref{rem:threshold}). The central limit theorem is applied only at the abscissa $x_0=\tfrac12+\delta_T$, which depends on $T$ but not on $t$. The passage from $x_0$ to the moving point $x_c$ is deterministic, at a cost $o(\sigma_T)$. For $f=\ell$ the divergence is extremely slow: at $t=10^{100}$ the standard deviation of $S(t;x_c)$ is only about $0.38$ (see the discussion after Conjecture~\ref{conj:gramlog}).

\begin{proof}[Proof of Theorem~\ref{thm:gramfail}]
Let $T$ be large, and put
\[
L=\log T,\quad\lambda=\log\log T,\quad L_f=L_f(T),\quad\delta=\delta_T,\quad x_0=\tfrac12+\delta,\quad\Psi=\Psi_T,
\]
\[
\mu=L_f/\sqrt\lambda,\quad B=\mu^{1/4},\quad W=B\sqrt\Psi/L,\quad h=4\pi/L .
\]
Then $\mu\to\infty$ and $B\to\infty$. Let $I=[x_0-2W,\,x_0+2W]$, $Q(t)=I\times[t-h,t+h]$ and $F(t)=\max_{Q(t)}|\zeta'/\zeta|$.

First, the formula for $\Psi$. We have $\sum_pp^{-1-v}=\log(1/v)+O(1)$ uniformly for $0<v\le1$, and $\sum_{p\ge T}p^{-1-2\delta}\le T^{-2\delta}/(2\delta L)=e^{-2L_f+O(1)}/(2L_f+O(1))\to0$. Since $2\delta=2L_f/L+O(1/L)$ and $L_f\to\infty$, this gives $\Psi=\tfrac12\log(L/L_f)+O(1)$. In particular $\Psi\le\tfrac12\lambda+O(1)$, so $B\sqrt\Psi\le\mu^{1/4}\sqrt\lambda=\mu^{-3/4}L_f$ for large $T$. Hence $2W\le\delta/2$. Note also that $\delta\ge4h$. Let $\mathcal Z$ be the set of zeros $\rho=\beta+i\gamma$ of $\zeta$ with $\beta>\tfrac12+\delta/4$.

We use the identities \eqref{eq:hadU} of Lemma~\ref{lem:hadamard}, and Selberg's density theorem \cite[Thm.~1]{Selberg46}: uniformly for $\tfrac12\le\sigma\le1$, the number of zeros with $\beta>\sigma$ and $0<\gamma\le3T$ is $\ll T^{1-(\sigma-1/2)/4}\log T$.

\emph{Step 1: typical behaviour.} Let $E_T$ be the set of $t\in[T,2T]$ such that
\begin{enumerate}
\item[(a)] $m_T(t)\le\tfrac23L_f$, with $m_T$ as in Lemma~\ref{lem:max};
\item[(b)] $|\zeta'/\zeta(x_0+it)|\le L/B^2$;
\item[(c)] $|\log|\zeta(x_0+it)||\le\tfrac14B\sqrt\Psi$;
\item[(d)] no $\rho\in\mathcal Z$ has $|\gamma-t|\le\delta$;
\item[(e)] $\widetilde N(t):=4\sum_{\rho\in\mathcal Z,\,|\gamma-t|>\delta}(\beta-\tfrac12)(\gamma-t)^{-2}\le L/8$;
\item[(f)] $\widetilde S(t):=4\sum_{\rho\in\mathcal Z,\,|\gamma-t|>\delta}(\gamma-t)^{-2}\le L^2/(B^3\sqrt\Psi)$.
\end{enumerate}
We claim that $\operatorname{meas}([T,2T]\setminus E_T)=o(T)$.

For (a), Lemma~\ref{lem:max} with $A=\tfrac23L_f$ bounds the exceptional set by $O(T\lambda/L_f^2)=O(T/\mu^2)$.

For (d)--(f), let $\mathcal Z'$ be the set of $\rho\in\mathcal Z$ with $0<\gamma\le3T$. Since $T^{-\delta/16}=e^{-L_f/16+O(1)}$, Selberg's theorem gives $\#\mathcal Z'\ll TLe^{-L_f/16}$ and
\[
\begin{split}
\sum_{\rho\in\mathcal Z'}(\beta-\tfrac12)&=\frac\delta4\,\#\mathcal Z'+\int_{1/2+\delta/4}^{1}\#\{\rho\in\mathcal Z':\beta>\sigma\}\,d\sigma\\
&\ll TL_fe^{-L_f/16}+TL\int_{\delta/4}^{\infty}e^{-vL/4}\,dv\ll T(L_f+1)e^{-L_f/16}.
\end{split}
\]
The zeros of $\mathcal Z$ outside $\mathcal Z'$ are at distance $\ge T$ from $[T,2T]$ and contribute $O(T^{-1}\log T)$ to $\widetilde N$ and $\widetilde S$. The set where (d) fails has measure at most $2\delta\,\#\mathcal Z'\ll TL_fe^{-L_f/16}$. Since $\int_{|v|>\delta}v^{-2}dv=2/\delta$, we have $\int_T^{2T}\widetilde N\ll\delta^{-1}\sum_{\mathcal Z'}(\beta-\tfrac12)\ll TLe^{-L_f/16}$ and $\int_T^{2T}\widetilde S\ll\delta^{-1}\#\mathcal Z'\ll TL^2e^{-L_f/16}/L_f$. So (e) and (f) fail on sets of measure $O(Te^{-L_f/16})$ and $O(TB^3\sqrt\Psi e^{-L_f/16}/L_f)=O(T\mu^{-1/4}e^{-L_f/16})$.

For (b) we use Lester's theorem \cite[Thm.~2]{LesterLD}, which is unconditional. Let $\psi(T)=(2\sigma-1)\log T$, $V=\tfrac12\sum_{n\ge2}\Lambda(n)^2n^{-2\sigma}$ and $\Omega=e^{-10}\min\bigl(V^{3/2},(\psi/\log\psi)^{1/2}\bigr)$. If $\psi(T)\to\infty$, $\psi(T)=o(\log T)$ and $r\Omega\ge1$, then
\[
\operatorname{meas}\Bigl\{t\in(0,T):\Bigl|\frac{\zeta'}{\zeta}(\sigma+it)\Bigr|\le r\sqrt V\Bigr\}=T\bigl(1-e^{-r^2/2}\bigr)+O\Bigl(T\,\frac{r^2+r}{\Omega}\Bigr).
\]
Take $\sigma=x_0$, so that $\psi(T)=2\delta L=2L_f+O(1)$, which tends to $\infty$ and is $o(\log T)$, and the same holds for $\psi(2T)$. We have $\Lambda(n)^2\le\Lambda(n)\log n$, and $\sum_n\Lambda(n)\log n\,n^{-s}=(\zeta'/\zeta)'(s)$, which is $(s-1)^{-2}+O(1)$ near $s=1$. Hence
\[
\frac{1}{2e^{2}}\sum_{p\le e^{1/\delta}}\frac{(\log p)^2}{p}\ \le\ V\ \le\ \frac12\Bigl(\frac{\zeta'}{\zeta}\Bigr)'(1+2\delta)=\frac{1+o(1)}{8\delta^2},
\]
and the left-hand side is $\gg\delta^{-2}$. Hence $\Omega\gg(L_f/\log L_f)^{1/2}\to\infty$, and $\sqrt V\le L/(2L_f)$ for large $T$. Let $r=\min(L_f/B^2,\Omega^{1/3})$. Then $r\to\infty$, $r^2=o(\Omega)$, and $r\sqrt V\le L/B^2$. Applying the theorem with $T$ and with $2T$, we see that (b) fails on a set of measure at most $2Te^{-r^2/2}+O(Tr^2/\Omega)=o(T)$.

For (c) we use the Selberg--Tsang central limit theorem near the critical line \cite{Tsang}, in the form stated in \cite{LesterLD}; see also \cite{HaLee} for sharper estimates on the lines $\sigma=\tfrac12+(\log T)^{-\theta}$. Let $0\le2\sigma-1=o(1)$ and $\Psi=\tfrac12\sum_{p<T}p^{-2\sigma}$. Then, uniformly for $a<b$,
\[
\operatorname{meas}\bigl\{t\in[0,T]:\ \Psi^{-1/2}\arg\zeta(\sigma+it)\in[a,b]\bigr\}=\frac{T}{\sqrt{2\pi}}\int_a^be^{-x^2/2}\,dx+O\Bigl(T\,\frac{\log\Psi}{\sqrt\Psi}\Bigr),
\]
and the same holds for $\log|\zeta(\sigma+it)|$ with error $O(T(\log\Psi)^2/\sqrt\Psi)$. We apply this with $\sigma=x_0$, so that $2\sigma-1=2\delta=o(1)$ and $\Psi=\Psi_T$, and with $T$ and $2T$. Since $\Psi\to\infty$ and $B\to\infty$, (c) fails on a set of measure $o(T)$. We shall also use the statement for $\arg\zeta$ in Step~4. Note that $\Psi_{2T}-\Psi_T=O(1/\log T)$ for the sums at the fixed abscissa $x_0$, and that $\arg\zeta(x_0+it)=\pi S(t;x_0)$ by definition.

\emph{Step 2: off-line zeros and a derivative bound.} Let $t\in E_T$ and $z=x+iv$ with $x\ge\tfrac12+\delta/4$ and $|v-t|\le h$. Every zero with $\beta>x$ lies in $\mathcal Z$, and by (d) it has $|\gamma-v|\ge|\gamma-t|-h\ge\tfrac12|\gamma-t|$. Hence, by (e),
\begin{equation}\label{eq:negpart}
\mathcal N(z):=\sum_{\beta>x}\frac{\beta-x}{|z-\rho|^2}\le\sum_{\rho\in\mathcal Z,\,|\gamma-t|>\delta}\frac{\beta-\frac12}{(\gamma-v)^2}\le\widetilde N(t)\le\frac L8 .
\end{equation}
In particular, by (d), $\zeta$ has no zeros $x+iv$ with $x>\tfrac12+\delta/4$ and $|v-t|\le h$. This region contains $Q(t)$, so $F(t)<\infty$. By \eqref{eq:hadU}, $\Re\frac{\zeta'}{\zeta}(z)\ge-\frac12\log\frac{|v|}{2\pi}-\frac L8-O(\frac1T)$.

We claim that $F(t)\le\eta_TL/B$, where $\eta_T\to0$ does not depend on $t$. Let $z\in Q(t)$, so $x-\tfrac12\ge\delta-2W\ge\delta/2$. In the second sum in \eqref{eq:hadU}, the zeros with $\beta>x-\delta/4$ lie in $\mathcal Z$ and contribute at most $\widetilde S(t)$, by (d). The zeros with $\beta\le x-\delta/4$ satisfy $|z-\rho|^{-2}\le(4/\delta)(x-\beta)|z-\rho|^{-2}$. By the first identity in \eqref{eq:hadU} and \eqref{eq:negpart}, their contribution is at most
\[
\frac4\delta\sum_{\beta<x}\frac{x-\beta}{|z-\rho|^2}=\frac4\delta\Bigl(\Re\frac{\zeta'}{\zeta}(z)+\frac12\log\frac{v}{2\pi}+\mathcal N(z)+O\Bigl(\frac1T\Bigr)\Bigr)\le\frac4\delta\bigl(F(t)+L\bigr).
\]
So $|(\zeta'/\zeta)'(z)|\le(4/\delta)(F(t)+L)+\widetilde S(t)+1$. The segment from $x_0+it$ to $z$ lies in $Q(t)$ and has length at most $2W+h$. With (b) this gives
\[
F(t)\le\frac{L}{B^2}+(2W+h)\Bigl(\frac4\delta\bigl(F(t)+L\bigr)+\widetilde S(t)+1\Bigr).
\]
Here $\eta_1:=(2W+h)\cdot4/\delta=(8B\sqrt\Psi+16\pi)/(L_f+O(1))\le9\mu^{-3/4}$ for large $T$, and $(2W+h)(\widetilde S(t)+1)\le3L/B^2$ by (f). Hence $F(t)\le(1-\eta_1)^{-1}(4L/B^2+\eta_1L)\le5L/B^2+2\eta_1L$. Since $\eta_1B\le9\mu^{-1/2}$, the claim holds with $\eta_T=5/B+18\mu^{-1/2}$. Consequently
\begin{equation}\label{eq:Fbound}
hF(t)\le4\pi\eta_T/B=o(1),\qquad 2WF(t)\le2\eta_T\sqrt\Psi=o(\sqrt\Psi).
\end{equation}

\emph{Step 3: good separators.} Let $t\in E_T$, and let $\tau=\tau_m$ be the separator nearest to $t$, so that $|\tau-t|\le h$. Write $L_\tau=\log(\tau/2\pi)$, and let $R$ and $\mathcal R$ be as in Lemma~\ref{lem:ident}. By (S2) there is a constant $C_1'$ with $|\log|f(\sigma'+i\tau)|-L_f|\le C_1'$ for $-1\le\sigma'\le2$ and $T-1\le\tau\le2T+1$.
\begin{enumerate}
\item[(i)] For $\tfrac12\le\sigma\le2$, (a) gives $|\zeta(\sigma+i\tau)|\le e^{2L_f/3}<\tfrac13|f(\sigma+i\tau)|$ for large $T$. So \eqref{eq:G1} holds.
\item[(ii)] $|\log\zeta(x_0+i\tau)-\log\zeta(x_0+it)|\le hF(t)=o(1)$. So by (c), $|\log|\zeta(x_0+i\tau)||\le\tfrac12B\sqrt\Psi$. Since there are no zeros in $[x_0,\infty)\times[t-h,t+h]$, also $|S(\tau;x_0)-S(t;x_0)|=o(1)$.
\item[(iii)] Single crossing. We verify the hypothesis of Lemma~\ref{lem:ident}(b) with $\eta=\delta/4$. For $\tfrac12<x\le\tfrac12+\delta/4$, Lemma~\ref{lem:ident}(a) and (a) give $\mathcal R(x)\le\tfrac14\delta L_\tau+\tfrac23L_f-L_f+C_1'+o(1)=-\tfrac1{12}L_f+O(1)<0$. For $\tfrac12+\delta/4<x\le2$, $\zeta(x+i\tau)\ne0$ and, by the bound after \eqref{eq:negpart}, $\Re\,\zeta'/\zeta(x+i\tau)\ge-\tfrac12L_\tau-\tfrac L8-O(1/T)\ge-\tfrac23L_\tau$. So Lemma~\ref{lem:ident}(b) applies. Thus $\mathcal R$ is increasing on $[\tfrac12+\delta/4,2]$ and has a unique zero $x_c$, and $k_m$ is the integer nearest to $-\tfrac12S(\tau;x_c)+e_m$ with $|e_m|\le\varepsilon_f(\tau)$. Next, $\delta L_\tau=L_f+O(1)$, so $\delta L_\tau-\log|f(1-x+i\tau)|=O(1)$ for $x\in I$, and by Lemma~\ref{lem:ident}(a)
\[
\mathcal R(x_0\pm2W)=\pm2WL_\tau+\log|\zeta(x_0\pm2W+i\tau)|+O(1).
\]
For $x\in I$ we have $|\log|\zeta(x+i\tau)|-\log|\zeta(x_0+i\tau)||\le2WF(t)=o(\sqrt\Psi)$, and $2WL_\tau\ge2B\sqrt\Psi(1-2/L)$. Hence $\pm\mathcal R(x_0\pm2W)\ge2B\sqrt\Psi(1-2/L)-\tfrac12B\sqrt\Psi-o(\sqrt\Psi)-O(1)>0$. As $I\subset[\tfrac12+\delta/4,2]$, it follows that $x_c\in I$.
\item[(iv)] By Step~2, $\zeta$ has no zeros on the half-line $\{x+i\tau:\ x\ge x_0-2W\}$, so $|S(\tau;x_c)-S(\tau;x_0)|\le2WF(t)/\pi=o(\sqrt\Psi)$.
\item[(v)] By (iii) and \eqref{eq:SHk}, $|S_f(\tau_m)-\tfrac12S(\tau;x_c)|\le|C_f|+\tfrac32$. With (ii) and (iv),
\[
\bigl|S_f(\tau_m)-\tfrac12S(t;x_0)\bigr|\le|C_f|+2+o(\sqrt\Psi)=o(\sigma_T).
\]
\end{enumerate}

\emph{Step 4: counting.} Consecutive separators in $[T-1,2T+1]$ are at most $g=2\pi/(L-2)$ apart. So if each $t\in E_T$ is assigned its nearest separator, each separator is assigned to a set of measure at most $g$. Note that $\tfrac12S(t;x_0)/\sigma_T=\arg\zeta(x_0+it)/\sqrt{\Psi_T}$. Let $-\infty\le a<b\le\infty$ and $\eta_0>0$. By Step~3(v), every $t\in E_T$ with $\arg\zeta(x_0+it)/\sqrt{\Psi_T}\in[a+\eta_0,b-\eta_0]$ has a nearest separator $\tau_m$ with $S_f(\tau_m)/\sigma_T\in[a,b]$, once $T$ is large. By Step~1 and the Selberg--Tsang theorem (applied with $T$ and $2T$), writing $N_T(a,b)$ for the number of $m$ with $T-h\le\tau_m\le2T+h$ and $S_f(\tau_m)\in[a\sigma_T,b\sigma_T]$, we get
\begin{align*}
N_T(a,b)&\ge\frac1g\operatorname{meas}\Bigl\{t\in E_T:\ \frac{\arg\zeta(x_0+it)}{\sqrt{\Psi_T}}\in[a+\eta_0,b-\eta_0]\Bigr\}\\
&\ge\Bigl(\frac1{\sqrt{2\pi}}\int_{a+\eta_0}^{b-\eta_0}e^{-x^2/2}\,dx-o(1)\Bigr)\frac{TL}{2\pi}.
\end{align*}
The number of all $m$ with $T<\tau_m\le2T$ is $(\theta(2T)-\theta(T))/\pi+O(1)=(1+o(1))TL/2\pi$, and the end ranges contain $O(1)$ separators. Applying the lower bound to $[a,b]$, to $(-\infty,a-\eta_0]$ and to $[b+\eta_0,\infty)$, and letting $\eta_0\to0$, gives the stated limit. The remaining statements follow since $\int_{-\epsilon}^{\epsilon}e^{-x^2/2}dx\to0$ as $\epsilon\to0$ and $\sigma_T\to\infty$, together with \eqref{eq:SHk} and Corollary~\ref{cor:block}.
\end{proof}

\begin{remark}[The threshold]\label{rem:threshold}
The hypothesis $L_f(T)/\sqrt{\log\log T}\to\infty$ is used three times: in Step~1(a), where Lemma~\ref{lem:max} is applied at height $\tfrac23L_f(T)$; to keep the window $I$, of width $4B\sqrt{\Psi_T}/\log T$, inside $[\tfrac12+\delta_T/4,2]$; and in the derivative bound of Step~2, which needs $B^2\sqrt{\Psi_T}=o(L_f(T))$. It marks the point where the balance-point description stops working. If $L_f(T)\le C\sqrt{\log\log T}$, then by the Selberg--Tsang theorem there is $c(C)>0$ such that $|\zeta(\sigma_T'+it)|>|f(\sigma_T'+it)|$, where $\sigma_T'=\tfrac12+1/\log T$, on a subset of $[T,2T]$ of measure at least $c(C)T$; recall that $\log|f(\sigma_T'+it)|=L_f(T)+O(1)$ there. At such heights $|\zeta/f|$ exceeds $1$ at a point to the right of the critical line, so the description of the counter by Proposition~\ref{prop:ident} and Lemma~\ref{lem:ident} is not available. We do not know whether the conclusions of Theorems~\ref{thm:right} and~\ref{thm:gramfail} fail there; see Section~\ref{sec:open}.
\end{remark}

\begin{corollary}[Neither]\label{cor:neither}
There is an absolute constant $c>0$ such that, for all large $T$,
\[
\#\bigl\{m:\ T<\tau_m\le2T,\ |N_H(\tau_m)-N(\tau_m)|\ge c\sqrt{\log\log T}\bigr\}\ge c\,T\log T .
\]
In particular $N_H-N$ is unbounded, and so is $S_H$.
\end{corollary}
\begin{proof}
Recall $N_H-N=S_H-1-S$, and let $\sigma_T$ be as in Theorem~\ref{thm:gramfail}, with $f=\ell$ and $\delta_T=\log\log T/\log T$. The Selberg--Tsang central limit theorem with $\sigma=\tfrac12$ (stated in Step~1 of the proof of Theorem~\ref{thm:gramfail}; here $\Psi=\tfrac12\log\log T+O(1)$) shows that $S(t)/(2\sigma_T)$ is asymptotically standard Gaussian for $t$ uniform in $[T,2T]$, since $4\sigma_T^2\sim\frac1{2\pi^2}\log\log T$; this is Selberg's theorem \cite{Selberg46}. We transfer this to separators. Let $Z_m$ be the number of zeros of $\zeta$, with multiplicity, with $\tau_{m-1}<\gamma\le\tau_m$. For $t\in[\tau_{m-1},\tau_m]$ we have $|S(t)-S(\tau_m)|\le Z_m+1$. Since $\sum_{T<\tau_m\le2T}Z_m\le(1+o(1))\frac{T}{2\pi}\log T$, at most a proportion $1/K+o(1)$ of these $m$ have $Z_m>K$. The intervals $[\tau_{m-1},\tau_m]$ with $Z_m>K$ have total length $O(T/K)$. Hence, arguing as in Step~4 of the proof of Theorem~\ref{thm:gramfail} and letting $K\to\infty$, the proportion of $m$ with $T<\tau_m\le2T$ and $|S(\tau_m)|>3\sigma_T$ is at least $\mathbb P(|Z|>\tfrac32)-o(1)>0.13$, where $Z$ is standard Gaussian. By Theorem~\ref{thm:gramfail}, the proportion with $|S_H(\tau_m)|>2.9\,\sigma_T$ tends to $\mathbb P(|Z|>2.9)<0.004$. If $|N_H(\tau_m)-N(\tau_m)|<c\sqrt{\log\log T}$, then $|S(\tau_m)-S_H(\tau_m)|<0.1\,\sigma_T$ for large $T$, provided $c\le1/100$. So the proportion of $m$ with $|N_H(\tau_m)-N(\tau_m)|\ge c\sqrt{\log\log T}$ is at least $0.12$ for large $T$, and the claim follows with $c=1/100$, since the number of $m$ with $T<\tau_m\le2T$ is $\sim\frac{T}{2\pi}\log T$.
\end{proof}

\begin{remark}[What $S_H$ is]\label{rem:SH}
By Lemma~\ref{lem:ident} and Step~3 of the proof of Theorem~\ref{thm:gramfail}, at almost all separators,
\[
S_H(\tau_m)=\tfrac12S\bigl(\tau_m;\ \tfrac12+\delta_m\bigr)+O(1),\qquad\delta_m=\frac{\log\log\tau_m-\log|\zeta(\tfrac12+\delta_m+i\tau_m)|}{\log(\tau_m/2\pi)}+O\Bigl(\frac1{\log^2\tau_m}\Bigr).
\]
On RH, Lemma~\ref{lem:ident}(c) gives this at every separator where $|\zeta|\le\frac13|\ell|$ on $[\tfrac12,2]$. So the counting error for $\zeta=\log(-s)$ is half of the argument function of $\zeta$, evaluated a distance $\delta_m$ to the right of the critical line. Since $S(t;\tfrac12+\delta)$ is a smoothing of $S(t)$, \eqref{eq:twoerr} suggests $N_H-N\approx-\tfrac12S$. Corollary~\ref{cor:neither} confirms that $N_H-N$ is indeed of size $\sqrt{\log\log T}$ at a positive proportion of separators.
\end{remark}

\subsection{The drift error}
In the rest of this section $f=\ell$.

\begin{proposition}[The drift error]\label{prop:Dformula}
Let
\[
D(T)=\sum_{0<\gamma\le T}(\beta-\tfrac12)+\frac1{2\pi}\int_2^T\log|1-\log(-it)|\,dt,\qquad Q(T)=\sum_{n\ge2}\frac{\Lambda(n)}{n^2\log^2n}\sin(T\log n).
\]
Then there is a constant $C$ such that, as $m\to\infty$,
\[
2\pi D(\tau_m)=3\pi k_m+\int_{-1}^{2}\bigl(\arg H(\sigma+i\tau_m)-\arg H(2+i\tau_m)\bigr)\,d\sigma+Q(\tau_m)+C+o(1),
\]
where $\arg H$ is continued horizontally from $2+i\tau_m$.
\end{proposition}
\begin{proof}
Apply Littlewood's lemma on $[-1,2]\times[T_0,\tau_m]$, keeping every term.
\begin{itemize}
\item Stirling's formula gives $\log|\chi(-1+it)|=\tfrac32\log(t/2\pi)+O(t^{-2})$, with no $1/t$ term, and $\theta'(t)=\tfrac12\log(t/2\pi)+O(t^{-2})$.
\item Exactly, $\int_{T_0}^{T}\log|\zeta(2-it)|\,dt=Q(T)-Q(T_0)$, and $\int\log|1-\ell/\zeta|(-1+it)\,dt$ converges.
\item On $\sigma=2$, the $u$-series of Lemma~\ref{lem:right} contributes a constant plus $O(1/\log T)$. Replacing $\ell(2+it)$ by $\log(-it)$ changes $\log|1-\ell|$ by $\Re\bigl(2i/(t(1-\log(-it)))\bigr)+O(t^{-2})=O(1/(t\log^2t))$, which is integrable.
\item The value $\arg H(2+i\tau_m)$ is a constant plus $o(1)$.
\item Finally, substitute $N_H(\tau_m)=m+C_H-k_m$ from Proposition~\ref{prop:counter}.
\end{itemize}
\end{proof}

\begin{remark}\label{rem:drift}
Proposition~\ref{prop:Dformula} expresses the drift error at the separators through the data that govern Gram's law. Besides the term $3\pi k_m$, the integrand $\arg H(\sigma+i\tau_m)-\arg H(2+i\tau_m)$ equals $0$ at $\sigma=2$ and $\alpha_m-2\pi k_m+o(1)$ at $\sigma=-1$, by Proposition~\ref{prop:counter}. In the computed range the counter is constant (Section~\ref{sec:numerics}, item~2), and $D$ stays in a short interval at the separators (item~4). We do not study the fluctuations of $D$ here.
\end{remark}

\subsection{Conjectures}

\begin{conjecture}\label{conj:gramlog}
For $\zeta(s)=\log(-s)$:
\begin{enumerate}
\item[(a)] The proportion of $m\le M$ with $n_m\ne1$ tends to $0$.
\item[(b)] Near height $T$, failures occur at the crossings of odd integers by $t\mapsto S(t;x)$, where $x=1-\bar\beta(T)$ and $\bar\beta(T)=\tfrac12-\log|1-\log(-iT)|/\log(T/2\pi)$. By Rice's formula, the expected proportion of failing intervals is
\[
p(T)\approx\frac{2\,s'}{s\log(T/2\pi)}\sum_{j\in\mathbb Z}e^{-(2j+1)^2/2s^2},
\]
where $s^2=\frac{1}{2\pi^2}\sum_{p,k}k^{-2}p^{-2kx}$ and $s'^2=\frac1{2\pi^2}\sum_n\Lambda(n)^2n^{-2x}$. As $T\to\infty$, $p(T)\sim(2\sqrt\pi\,\log\log T)^{-1}$.
\end{enumerate}
\end{conjecture}

Both Theorem~\ref{thm:gramfail} and the conjecture take effect extremely slowly, and nothing in the proof of Theorem~\ref{thm:gramfail} is effective at accessible heights. The variance $\Psi$ of $\pi S(t;x_0)$ is $\asymp\log\log T$, but its actual size is small: with $x_0=\tfrac12+\log\log T/\log T$ one finds $\Psi=\tfrac12\sum_{p<T}p^{-2x_0}\approx0.44$, $0.58$, $0.84$ and $1.40$ at $T=10^4$, $10^8$, $10^{20}$ and $10^{100}$, so the standard deviation of $S(t;x_0)$ is only about $0.21$, $0.24$, $0.29$ and $0.38$. Moreover, a failure requires a \emph{crossing} of an odd integer by $S(\cdot\,;x_c)$, not merely a value beyond one. The formula in (b) gives $p(T)\approx5\cdot10^{-6}$ at $T=10^4$ and $1.4\cdot10^{-4}$ at $T=10^8$, that is, $0.01$ and $0.27$ expected failures per window of $2000$ intervals (Table~\ref{tab:gram}). The rate is about $4\cdot10^{-3}$ at $T=10^{50}$ and $7\cdot10^{-3}$ at $T=10^{100}$. It peaks at about $1.8\cdot10^{-2}$ for $T$ between $10^{10^4}$ and $10^{10^5}$, and then decays like $1/\log\log T$. The model also ignores the selection effect described in Remark~\ref{rem:selection}, which makes actual failures rarer still at these heights.

\section{Numerical evidence}\label{sec:numerics}

\subsection*{Methods}
All computations used Arb~\cite{Arb} through python-flint, together with mpmath. We distinguish two levels of rigour.

\emph{Interval arithmetic.} The number of zeros of $H$ in $[-1.5,3]\times[1,T]$, for $T=100,200,\dots,2000$, was computed by a ball-arithmetic argument principle at $64$-bit precision. The contour is cut into segments. On each segment $[z_0,z_1]$, Arb evaluates $H$ on the disc with diameter $[z_0,z_1]$ and returns a rectangle enclosing $H$ of that disc. If the enclosing rectangle avoids $0$, it lies in an open half-plane through $0$. The continuous change of $\arg H$ along the segment is then less than $\pi$ in absolute value, so it equals the principal difference of the arguments at the endpoints. If the rectangle meets $0$, the segment is bisected. Summing over segments gives the winding number exactly; the sums agreed with integers to within $10^{-6}$.

The same procedure was applied separately to each of the $2000$ shifted Gram rectangles $[-1,2]\times[\tau_{m-1},\tau_m]$ starting at height $10^4$. The endpoints $\tau_m$ were computed with mpmath to $30$ digits. Every rectangle has winding number exactly $1$.

\emph{Floating point.} Individual zeros were located by Newton's method in floating point, with residual $|H|<10^{-8}$. The interval-arithmetic counts at each $T\in100\mathbb Z$ equal the number of distinct zeros located, so no zero in $[-1.5,3]\times[1,2000]$ was missed. The per-interval counts at heights $10^5,\dots,10^8$ (Table~\ref{tab:gram}) used the same contour decomposition with double-precision samples and adaptive refinement, but without ball enclosures, so they are numerical rather than rigorous.

\emph{Counting convention.} The theorems count the solutions with $-1<\sigma<2$ and $\gamma>0$ (Section~\ref{sec:slow}), which for large $\gamma$ are all of them (Lemma~\ref{lem:strip}). Only finitely many lie outside the half-strip $\mathcal S=[-1.5,3]\times[1,\infty)$, and the statistics below are for the solutions in $\mathcal S$. In this section $N_H(T)$ counts the solutions in $\mathcal S$ with $\gamma\le T$, and $\rho_n=\beta_n+i\gamma_n$ is the $n$-th of them by height. A floating-point argument principle over $[-300,3]\times[10^{-6},200]$ and $[-40,-1.5]\times[200,2000]$ finds exactly four solutions with $\gamma>0$ outside $\mathcal S$ there:
\[
-15.7274+1.2145i,\qquad-12.3995+3.1478i,\qquad-8.0835+5.2817i,\qquad0.9295+0.3406i .
\]
There are none with $\sigma\ge3$, since there $|\zeta(s)|\le\zeta(3)<\tfrac\pi2<|\log(-s)|$. Counting these four as well, the number of solutions with $0<\gamma\le T$ is $N_H(T)+4$ for $T\ge5.29$; every offset $n-\theta(\gamma_n)/\pi$ below then increases by $4$, and every sum of $\beta-\tfrac12$ decreases by the constant $37.28$.

\begin{table}[htbp]
\centering
\caption{Window means of $\beta$ for zeros of $H$, compared with $\tfrac12-\int\log|1-\log(-it)|\,dt\big/\int\log(t/2\pi)\,dt$ over the window. For the first window the integrals start at the first ordinate, $11.34$.}\label{tab:drift}
\begin{tabular}{rrrrr}
\toprule
window & $n$ & mean $\beta$ & predicted & share with $\beta>\frac12$\\
\midrule
$(0,250]$ & 108 & 0.0075 & 0.0152 & 4.63\%\\
$(250,500]$ & 162 & 0.0989 & 0.0973 & 3.70\%\\
$(500,750]$ & 183 & 0.1218 & 0.1229 & 4.37\%\\
$(750,1000]$ & 196 & 0.1377 & 0.1374 & 4.59\%\\
$(1000,1250]$ & 206 & 0.1463 & 0.1474 & 3.88\%\\
$(1250,1500]$ & 215 & 0.1558 & 0.1548 & 3.72\%\\
$(1500,1750]$ & 221 & 0.1616 & 0.1607 & 3.62\%\\
$(1750,2000]$ & 226 & 0.1644 & 0.1656 & 3.10\%\\
\bottomrule
\end{tabular}
\end{table}

\begin{table}[htbp]
\centering
\caption{Zeros of $\zeta(s)-\log(-s)$ per shifted Gram interval (argument principle), $2000$ consecutive intervals per window. Expected failures per window are from Conjecture~\ref{conj:gramlog}(b).}\label{tab:gram}
\begin{tabular}{rrrl}
\toprule
start height & $1-\bar\beta$ & expected failures & observed\\
\midrule
$10^4$ & 0.788 & $0.009$ & all 2000 intervals contain exactly 1\\
$10^5$ & 0.744 & $0.035$ & all 2000 intervals contain exactly 1\\
$10^6$ & 0.714 & $0.087$ & all 2000 intervals contain exactly 1\\
$10^7$ & 0.691 & $0.166$ & all 2000 intervals contain exactly 1\\
$10^8$ & 0.673 & $0.273$ & all 2000 intervals contain exactly 1\\
\bottomrule
\end{tabular}
\end{table}

\begin{figure}[htbp]
\centering
\includegraphics[width=0.9\linewidth]{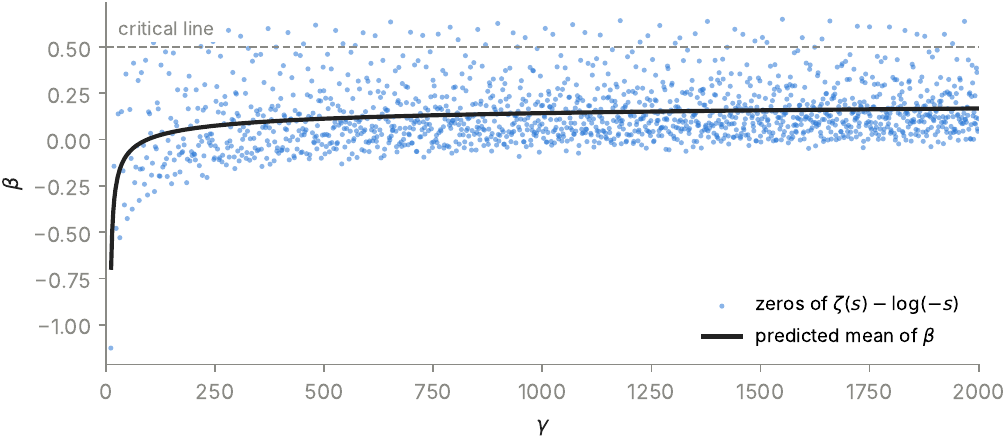}
\caption{Real parts of the $1517$ zeros of $\zeta(s)-\log(-s)$ in $[-1.5,3]\times[1,2000]$, with the predicted mean $\tfrac12-\log|1-\log(-it)|/\log(t/2\pi)$.}\label{fig:beta}
\end{figure}

\begin{figure}[htbp]
\centering
\includegraphics[width=0.9\linewidth]{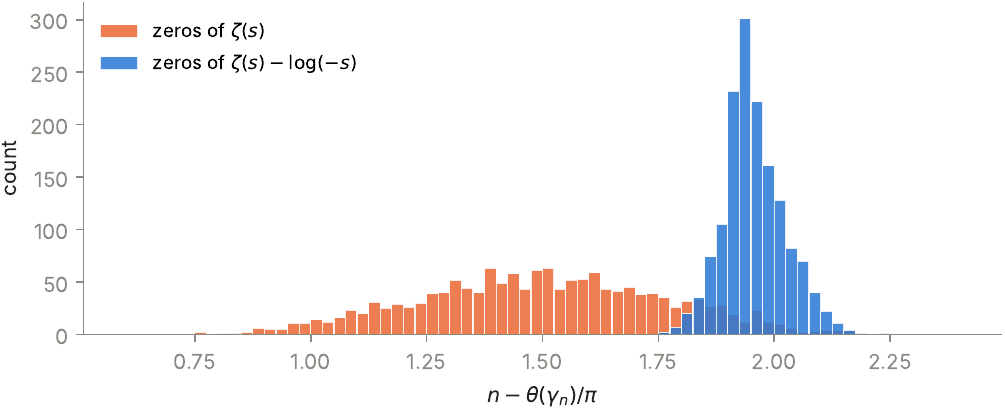}
\caption{$n-\theta(\gamma_n)/\pi$ for the zeros of $\zeta(s)-\log(-s)$ in $[-1.5,3]\times[1,2000]$ and for the zeros of $\zeta(s)$ with $\gamma\le2000$. The solutions lock to Gram points; the zeta zeros do not.}\label{fig:phase}
\end{figure}

\begin{figure}[htbp]
\centering
\includegraphics[width=0.9\linewidth]{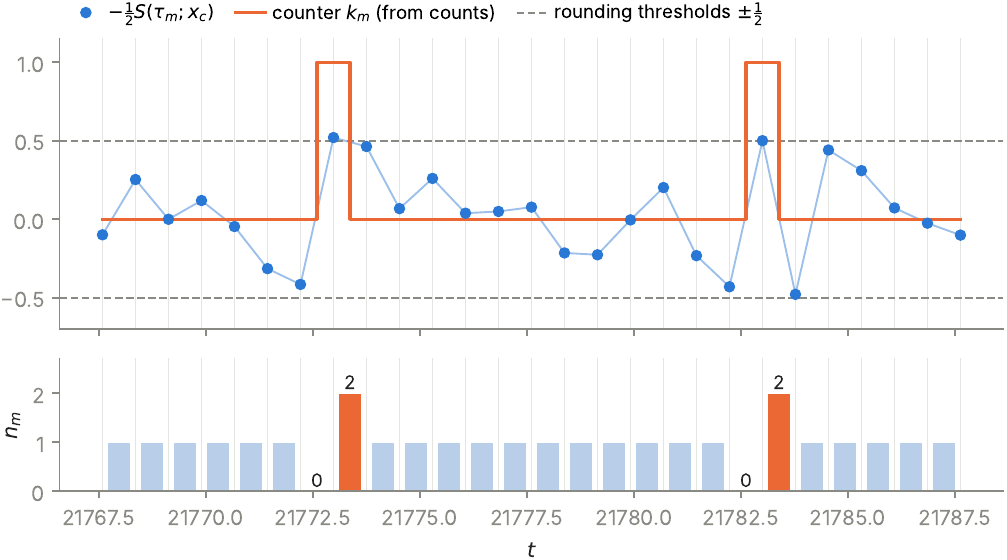}
\caption{Two phase slips for $f(s)=(-s/2\pi)^{0.05}$ near $t=21{,}775$. Top: $-\tfrac12S(\tau_m;x_c)$ at the separators (dots) and the counter $k_m$ computed from the solution counts (step). Bottom: $n_m$. When the phase passes $\tfrac12$, $k_m$ jumps to $1$, and $n_m=1+k_{m-1}-k_m$ gives an empty interval followed by a double one.}\label{fig:slip}
\end{figure}

\begin{enumerate}
\item \textbf{Counts.} $N_H(2000)=1517=N(2000)$. Moreover $N_H(T)-N(T)\in\{-1,0,1\}$ for all $T\le2000$. The ordinates interlace with the zeta zeros except at $8$ places.
\item \textbf{Gram's law.} $\operatorname{round}(\theta(\gamma_n)/\pi)=n-2$ for every $n\le1517$, that is, $n-\theta(\gamma_n)/\pi\in[1.77,2.19]$. For the zeta zeros the same statistic ranges over $[0.77,2.23]$ (Figure~\ref{fig:phase}). The residual in the phase identity \eqref{eq:phase}, multiplied by $\gamma$, is at most $1.32$.
\item \textbf{Higher windows.} Starting at $T=10^4,\dots,10^8$, we took $2000$ consecutive shifted Gram intervals and counted zeros of $H$ in each by the argument principle. Every interval contains exactly one zero (Table~\ref{tab:gram}). Conjecture~\ref{conj:gramlog}(b) predicts between $0.01$ and $0.27$ failures per window, so observing none is expected, and these data do not yet test the conjecture.
\item \textbf{Drift law.} The drift error $D(T)$ of Proposition~\ref{prop:Dformula}, summed over the solutions in $\mathcal S$, is a step function that jumps by $\beta-\tfrac12$ at each ordinate. At the $1495$ separators $\tau_m\in[60,1990]$ it lies in $[-0.46,0.12]$, with standard deviation $0.14$. Over all $T\le2000$ it lies in $[-0.75,0.88]$, while $\sum(\beta-\tfrac12)$ itself reaches $-555.9$. Window means of $\beta$ agree with the prediction to within $0.002$ beyond $t=250$ (Table~\ref{tab:drift} and Figure~\ref{fig:beta}).
\item \textbf{Crossover structure.} On each of $1500$ random horizontal lines with $1000\le t\le2000$, the function $|\zeta(s)|/|\log(-s)|$ crosses $1$ exactly once, as Lemma~\ref{lem:ident} predicts on RH. Tracking the crossover phase $\alpha_c(t)$ on a grid and solving $\alpha_c(t)\equiv2\theta(t)$ reproduces all $334$ zeros of $H$ in $[1000,1400]$, with ordinates accurate to $2.4\cdot10^{-5}$.
\item \textbf{Right of the line.} $59$ of the $1517$ zeros have $\beta>\tfrac12$, and $\sum_{\beta>1/2}(\beta-\tfrac12)=4.17$ at $T=2000$.
\item \textbf{Growing targets.} For $f(s)=s$, $(-s)^{3/4}$ and $(-s)^{1/4}$, all $400$ shifted Gram intervals near $T=10^3$ contain exactly one solution. The same holds for $1000$ intervals near $T=10^5$ when $f(s)=(-s)^{1/4}$. Theorem~\ref{thm:gram} predicts the offsets $n-\theta(\gamma_n)/\pi$ modulo $1$, as $\varphi_0/2\pi$. With $n$ counting the solutions in $[-3,3]\times(0.05,\gamma_n]$, the offsets have mean $2.812$ and range $[2.71,2.93]$ for $\kappa=\tfrac34$ (predicted $\equiv0.8125$), mean $2.926$ and range $[2.74,3.16]$ for $\kappa=\tfrac14$ (predicted $\equiv0.9375$), and mean $1.250$ and range $[1.17,1.32]$ for $f(s)=s$ (predicted $\equiv0.25$). For $\kappa=\tfrac14$, Theorem~\ref{thm:smallkappa} predicts failures only where $|\arg\zeta(\tfrac34+it)|>\pi$, which is very rare; none was seen. With the much smaller exponent in $f(s)=(-s/2\pi)^{0.05}$, $112$ of $3000$ consecutive shifted Gram intervals starting at $T=2\cdot10^4$ fail: $57$ contain no solution, $53$ contain two and $2$ contain three. Of the failures, $106$ form adjacent pairs with counts $0$ and $2$: a solution moves into a neighbouring interval and back (Figure~\ref{fig:slip}). For this target the crossing of Proposition~\ref{prop:ident} is not unique at $122$ of the $3001$ separators at this height; at the $140$ separators with a unique crossing among the $142$ in $[21681,21790]$, $k_m$ is the integer nearest to $-\tfrac12S(\tau_m;x_c)$, as predicted.
\end{enumerate}

\begin{remark}[The selection effect]\label{rem:selection}
At the zeros of $H$ with $1000<\gamma\le2000$, $\arg\zeta(1-\rho)$ has standard deviation $0.435$, against about $0.65$ at the same abscissae and random heights. The zeros are where the crossover phase $\alpha_c(t)$ meets $2\theta(t)$ modulo $2\pi$, so by the Kac--Rice formula they sample $\alpha_c$ conditionally on that event; and $\alpha_c$ depends on $2\theta$ through the reflected sum $\chi(s)^{-1}\sum_{n\le\sqrt{t/2\pi}}n^{-s}$ in the approximate functional equation for $\zeta(1-s)$, which carries the phase $e^{2i\theta(t)}$ and has relative size about $(\log t)^{-1/2}$ for $f=\ell$. Over $[1000,1400]$ the spread of $\alpha_c$ is $0.48$ when the Gram phase is within $\pi/4$ of $0$ and $0.71$ when it is within $\pi/4$ of $\pi$. The ratio of the two spreads is $0.678\pm0.016$ in $[10^3,2\cdot10^3]$ and $0.862\pm0.020$ in $[10^6,10^6+2000]$, so this effect, which makes failures rarer than Conjecture~\ref{conj:gramlog}(b) predicts, fades as $T\to\infty$.
\end{remark}

\section{Open problems}\label{sec:open}
\begin{enumerate}\setlength{\itemsep}{0pt}
\item Prove Conjectures~\ref{conj:boundary} (the case $\kappa=\kappa_1$) and~\ref{conj:gramlog}.
\item For $0<\kappa\le\tfrac12$, find the order of the number of failures up to height $T$.
\item If $\log|f(iT)|\ll\sqrt{\log\log T}$, does a positive proportion of the solutions lie to the right of the critical line, and does $S_f(\tau_m)$ still satisfy a central limit theorem?
\end{enumerate}


\clearpage
\begin{thebibliography}{99}
\bibitem{BinderPauliSaidak} T.~Binder, S.~Pauli and F.~Saidak, Zeros of high derivatives of the Riemann zeta function, \emph{Rocky Mountain J. Math.} 45 (2015), 903--928.
\bibitem{BohrJessen} H.~Bohr and B.~Jessen, \"Uber die Wertverteilung der Riemannschen Zetafunktion, \emph{Acta Math.} 54 (1930), 1--35; Zweite Mitteilung, \emph{Acta Math.} 58 (1932), 1--55.
\bibitem{CardwellYe} M.~Cardwell and Z.~Ye, Value distribution of $L$-functions with rational moving targets, \emph{Adv. Pure Math.} 3 (2013), 719--723.
\bibitem{Guo} C.~R.~Guo, The distribution of the logarithmic derivative of the Riemann zeta function, \emph{Proc. London Math. Soc.} (3) 72 (1996), 1--27.
\bibitem{HaLee} J.~Ha and Y.~Lee, The $a$-values of the Riemann zeta function near the critical line, arXiv:1711.08928.
\bibitem{HRS} W.~Heap, M.~Radziwi\l\l\ and K.~Soundararajan, Sharp upper bounds for fractional moments of the Riemann zeta function, \emph{Q. J. Math.} 70 (2019), no.~4, 1387--1396.
\bibitem{HS} W.~Heap and K.~Soundararajan, Lower bounds for moments of zeta and $L$-functions revisited, \emph{Mathematika} 68 (2022), 1--14.
\bibitem{Arb} F.~Johansson, Arb: efficient arbitrary-precision midpoint-radius interval arithmetic, \emph{IEEE Trans. Comput.} 66 (2017), 1281--1292.
\bibitem{KalpokasSteuding} J.~Kalpokas and J.~Steuding, On the value-distribution of the Riemann zeta-function on the critical line, arXiv:0907.1910.
\bibitem{Laurincikas} A.~Laurin\v{c}ikas, \emph{Limit Theorems for the Riemann Zeta-Function}, Kluwer, Dordrecht, 1996.
\bibitem{LesterLD} S.~J.~Lester, The distribution of the logarithmic derivative of the Riemann zeta-function, arXiv:1308.3597.
\bibitem{Lester} S.~J.~Lester, $a$-points of the Riemann zeta-function on the critical line, \emph{Int. Math. Res. Not. IMRN} 2015, doi:10.1093/imrn/rnt356.
\bibitem{Levinson} N.~Levinson, Almost all roots of $\zeta(s)=a$ are arbitrarily close to $\sigma=1/2$, \emph{Proc. Nat. Acad. Sci. USA} 72 (1975), 1322--1324.
\bibitem{LiSteuding} B.~Q.~Li and J.~Steuding, Fixed points of the Riemann zeta function and Dirichlet series, \emph{Monatsh. Math.} 198 (2022), 581--589.
\bibitem{LiSteudingSuzuki} B.~Q.~Li, J.~Steuding and Y.~Suzuki, An asymptotic estimate for the characteristic and number of fixed points of the Riemann zeta function, \emph{Nagoya Math. J.} 260 (2025), 648--665; doi:10.1017/nmj.2025.10066.
\bibitem{MSZ} Yu.~Matiyasevich, F.~Saidak and P.~Zvengrowski, Horizontal monotonicity of the modulus of the zeta function, $L$-functions, and related functions, \emph{Acta Arith.} 166 (2014), 189--200.
\bibitem{Selberg46} A.~Selberg, Contributions to the theory of the Riemann zeta-function, \emph{Arch. Math. Naturvid.} 48 (1946), 89--155.
\bibitem{SSS} A.~Sourmelidis, J.~Steuding and A.~I.~Suriajaya, The $a$-points of the Riemann zeta-function and the functional equation, arXiv:2204.13887.
\bibitem{Stein} E.~M.~Stein, \emph{Singular Integrals and Differentiability Properties of Functions}, Princeton University Press, Princeton, 1970.
\bibitem{SteudingSuriajaya} J.~Steuding and A.~I.~Suriajaya, Value-distribution of the Riemann zeta-function along its Julia lines, \emph{Comput. Methods Funct. Theory} 20 (2020), 389--401.
\bibitem{Titchmarsh} E.~C.~Titchmarsh, \emph{The Theory of the Riemann Zeta-Function}, 2nd ed., Oxford, 1986.
\bibitem{Trudgian} T.~Trudgian, On the success and failure of Gram's law and the Rosser rule, \emph{Acta Arith.} 148 (2011), 225--256.
\bibitem{Tsang} K.-M.~Tsang, \emph{The distribution of the values of the Riemann zeta-function}, PhD thesis, Princeton University, 1984.
\bibitem{WangWangYao} J.~Wang, Q.~Wang and X.~Yao, On uniqueness of functions in the extended Selberg class with moving targets, arXiv:2502.09091.
\bibitem{Yamanoi} K.~Yamanoi, The second main theorem for small functions and related problems, \emph{Acta Math.} 192 (2004), 225--294.
\end{thebibliography}
\end{document}